\documentclass[11pt,a4paper]{amsart}

\makeatletter
 \def\@textbottom{\vskip \z@ \@plus 100pt}
 \let\@texttop\relax
\makeatother

\usepackage{amscd,amssymb,latexsym,url,verbatim,graphicx,color}

\usepackage{cases,amsmath,amsthm}

\usepackage{tikz,tikz-cd}

\usepackage[dvips]{epsfig}

\usepackage{txfonts}

\title
{Quantitative aspects of acyclicity}

\author{Dmitry N. Kozlov and Roy Meshulam}

\address{Dmitry N. Kozlov \\ Department of Mathematics, University of
  Bremen, 28334 Bremen, Federal Republic of Germany}

\email{dfk@math.uni-bremen.de}

\address{Roy Meshulam \\ Department of Mathematics, Technion - Israel
  Institute of Technology, Haifa 32000, Israel}

\email{meshulam@math.technion.ac.il}

\newtheorem{theorem}{Theorem}[section]

\newtheorem{df}[theorem]{Definition}
\newtheorem{thm}[theorem]{Theorem} 
\newtheorem{prop}[theorem]{Proposition}
\newtheorem{lm}[theorem]{Lemma}
\newtheorem{claim}[theorem]{Claim}
\newtheorem{crl}[theorem]{Corollary}

\newtheorem{rem}[theorem]{Remark}

 \newcommand{\nin}{\noindent}
\newcommand{\pr}{\nin{\bf Proof.} }
\newcommand{\prn}[1]{\nin{\bf Proof of #1.}}

\newcommand{\cf}{{\mathcal F}}
\newcommand{\cg}{{\mathcal G}}

\newcommand{\es}{\emptyset}

\newcommand{\ra}{\rightarrow}
\newcommand{\sm}{\setminus}
\newcommand{\supp}{\text{\rm supp}\,}

\newcommand{\ti}{\tilde}
\newcommand{\wti}{\widetilde}

\newcommand{\dz}{{\mathbb Z}_2}
\newcommand{\sy} {\text{sys}}
\newcommand{\csy}{\text{csy}}
\newcommand{\rk}{\text{rk}}

\newcommand{\enp}{\qed\vskip5pt}
\newcommand{\beq}[0]{\begin{equation}}
\newcommand{\enq}[0]{\end{equation}}
\newcommand{\Rea}{\mathbb{R}}

\newcommand{\Int}{\mathbb{Z}}

\newcommand{\FF}{\mathbb{F}}
\newcommand{\eps}{{\epsilon}}
\newcommand{\dn}{{\smp{n-1}}}

\newcommand{\thh}{\widetilde{H}}
\newcommand{\xv}{X^{\vee}}
\newcommand{\yv}{Y^{\vee}}
\newcommand{\prob}{\text{Pr}}
\newcommand{\ueps}{\underline{\epsilon}}
\newcommand{\ddp}{\smp{p-1}}
\newcommand{\mkn}{\min\{k,n-2\}}
\newcommand{\tphi}{\tilde{\phi}}
\newcommand{\dist}{\text{dist}}
\newcommand{\diam}{\text{\rm diam}\,}
\newcommand{\Hom}{\text{\rm Hom}\,}
\newcommand{\hzero}{\widehat{0}}
\newcommand{\hone}{\widehat{1}}

\newcommand{\idd}{{\rm Id}}
\newcommand{\ccs}{{\mathbb S}}

\newcommand{\sn}[1]{\left|#1\right|}
\newcommand{\n}[1]{\left\|#1\right\|}
\newcommand{\csyn}[1]{\n{#1}_\csy}
\newcommand{\myrb}[1]{\left(#1\right)}
\newcommand{\mycb}[1]{\left\{#1\right\}}
\newcommand{\mybar}[1]{\left|#1\right|}

\newcommand{\xor}{\oplus}

\newcommand{\im}{\textrm{im}\,}

\newcommand{\ev}[2]{\langle#1,#2\rangle}

\newcommand{\mysubs}[1]{\subsection{#1} $\,$

}

\newcommand{\fy}{\varphi}
\renewcommand{\phi}{\varphi}

\newcommand{\editcut}{}

\newcommand{\smp}[1]{\Delta^{#1}}
\newcommand{\bou}[1]{\partial_{#1}}
\newcommand{\cob}[1]{d_{#1}}
\newcommand{\rbou}[1]{\ti\partial_{#1}}
\newcommand{\rcob}[1]{\ti d_{#1}}

\numberwithin{equation}{section}
\numberwithin{figure}{section}
\numberwithin{table}{section}

\begin{document}

\begin{abstract}
The Cheeger constant is a measure of the edge expansion of a graph, and as such plays a key role in combinatorics and theoretical computer science. In recent years there is an interest in $k$-dimensional versions of the Cheeger constant that likewise provide quantitative measure of cohomological acyclicity of a complex in dimension $k$. In this paper we study several aspects of the higher Cheeger constants. Our results include methods for bounding the cosystolic norm of $k$-cochains and the $k$-th Cheeger constants, with applications to the expansion of pseudomanifolds, Coxeter complexes and homogenous geometric lattices.
We revisit a theorem of Gromov on the expansion of a product of a complex with a simplex, and provide an elementary derivation of the expansion in a hypercube. We prove a lower bound on the maximal cosystole in a complex and an upper bound on the expansion of bounded degree complexes, and  give an essentially sharp estimate for the cosystolic norm of the Paley cochains.
Finally, we discuss a non-abelian version of the $1$-dimensional expansion of a simplex,
with an application to a question of Babson on bounded quotients of the fundamental group of a random $2$-complex.
\end{abstract}

\maketitle

\section{Introduction}

The Cheeger constant is a parameter that quantifies the edge expansion of a graph, and as such
plays a key role in combinatorics and theoretical computer science (see, e.g., \cite{HLW06,Lu2}).
The  $k$-dimensional version of the graphical Cheeger constant, called "coboundary expansion", came up independently in the work of Linial, Meshulam and Wallach \cite{LiM,MW} on homological connectivity of random complexes and in Gromov's remarkable work \cite{Gr} on the topological overlap property;
see also the paper by Dotterer and Kahle \cite{DK}.
Roughly speaking, the $k$-th coboundary expansion $h^k(X)$ of a polyhedral complex $X$ is a measure of the minimal distance of $X$ from a complex $Y$ that satisfies
$H^k(Y;\Int_2) \neq 0$. Likewise, $h_k(X)$ is a measure of the minimal distance of $X$ from a complex $Y$ that satisfies $H_k(Y;\Int_2) \neq 0$.
We proceed with the formal definitions.

\mysubs{Chains and cochains.}
Let $X$ be a~finite polyhedral complex. In this paper we shall mostly work with
$\dz$ coefficients. Accordingly, let $C_k(X)=C_k(X;\dz)$ denote the
space of $k$-chains of $X$ over $\dz$ and let $C^k(X)=C^k(X;\dz)$
denote the space of $\dz$-valued $k$-cochains of $X$.  Let
$\bou{k}:C_k(X) \rightarrow C_{k-1}(X)$ and $\cob{k}:C^k(X)\rightarrow
C^{k+1}(X)$ denote the usual $k$-boundary and $k$-coboundary
operators.

Let $\ev{}{}:C^k(X)\times C_k(X)\ra \dz$ denote the usual evaluation
of $k$-cochains on $k$-chains. Note that
\begin{equation}
\label{e:eval}
\ev{\phi}{\bou{k}c}=\ev{\cob{k-1}\phi}{c}
\end{equation}
whenever $c \in C_k(X)$ and $\phi \in C^{k-1}(X)$.

For each non-negative $k$, let $X^{(k)}$ denote the $k$-th skeleton of~$X$,
which is itself a~finite polyhedral complex. Furthermore, let $X(k)$ denote
the set of all $k$-cells of~$X$. The {\bf degree} of a face $\sigma \in X(k)$ is
$\deg_X(\sigma)=|\{\tau \in X(k+1): \sigma \subset \tau\}|$. In particular, if $G=(V,E)$ is a graph, then $\deg_G(v)$ is the usual degree of a vertex $v$.

When $X$ is a simplicial complex and $\{a_0,\ldots,a_k\} \in X$, we write $[a_0,\ldots,a_k]$ for the corresponding element in $C_k(X)$. As we are working
over $\dz$, the order of the $a_i$'s does not matter. We will use the convention that
$[a_0,\ldots,a_k]=0$ if $a_i=a_j$ for some $0 \leq i \neq j \leq k$.

For any set of cells $A$ and any cochain $\phi$, we let
$\phi_A$ denote the restriction cochain, i.e., the cochain that is
equal to $\phi$ on the set~$A$ and is equal to~$0$ otherwise.
Furthermore, we let $A^*$ denote the cochain which is equal to $1$ on
the set $A$ and is $0$ otherwise.

In our convention the set $X(-1)$ is empty, and accordingly $C_{-1}(X)=0$,
forcing $\bou{0}=0$ and $\cob{-1}=0$. This is the so-called non-reduced
setting. In the reduced setting we let $\wti X(-1)$ be the set containing
a~single element, denoted $\es_X$. This is the so-called empty simplex.
For convenience we also set $\wti X(k):=X(k)$, for $k\neq -1$.
Accordingly, the reduced boundary and coboundary operators coincide with
the non-reduced ones except for the following cases:
\begin{itemize}
\item $\rbou{0}(v)=\es_x$, for all $v\in X(0)$;
\item $\rcob{-1}(\es_X^*)=\sum_{v\in X(0)}v^*$.
\end{itemize}

For a subcomplex $Y\subset X$ let $C_k(X,Y)=C_k(X)/C_k(Y)$ denote
the space of relative $k$-chains, with its induced $k$-boundary map
$\partial_k$. Identifying $C^k(Y)$ with the subspace of $C^k(X)$
consisting of all $k$-cochains whose support is contained in $Y$, let
$C^k(X,Y)=C^k(X)/C^k(Y)$ denote the space of relative $k$-cochains,
with its induced $k$-coboundary map $d_k$.  Let
$B_k(X,Y)=\partial_{k+1}(C_{k+1}(X,Y))$ and
$B^k(X,Y)=d_{k-1}(C^{k-1}(X,Y))$ denote the spaces of relative
$k$-boundaries and relative $k$-coboundaries.

\mysubs{Homology expansion}
Let $X$ be a finite polyhedral complex $X$ and let $c \in C_k(X)$.
Write $c=\sum_{\sigma \in X(k)} a_\sigma \sigma$
where the $a_{\sigma}$'s are in $\dz$.

\begin{df} \label{df:snorm}
The {\bf norm} of $c$ is
\[
\|c\|:=|\supp c|=|\{\sigma \in X(k):a_{\sigma} \neq 0\}|.
\]
The {\bf systolic norm} of $c$ is
\[
\|c\|_{sys}:=\min\big\{\|c+\partial_{k+1}c'\|:c' \in C_{k+1}(X)\big\}.
\]
The chain $c$ is called a~{\bf $k$-systole} if $\|c\|_{sys}=\|c\|$.
A~{\bf systolic form} of $c$ is any
  $\tilde{c}=c+\partial_{k+1}c'$, such that
  $\|\ti{c}\|=\|c\|_\sy$.
\end{df}

\begin{df}
The~{\bf boundary expansion} of a~$k$-chain $c \in C_k(X) \setminus B_k(X)$ is
\[\|c\|_{\exp}:=\|\partial_k c\|/\|c\|_\sy.\]
The {\bf $k$-th homological Cheeger constant} of $X$ is
\[h_k(X):=\min_{c \in C_k(X) \setminus B_k(X)}\|c\|_{\exp}.\]
\end{df}

The definition of expansion can be extended to the relative case as follows.
Let $Y$ be a subcomplex of~$X$. A~relative chain in $C_k(X,Y)$ has a~unique representative
$c=\sum_{\sigma\in X(k)\setminus Y(k) }a_\sigma \sigma$.  The {\bf norm} of $c+C_k(Y)$
is then defined by $\|c+C_k(Y)\|:=\|c\|$. The notions of the relative systolic norm, relative expansion and
relative homological Cheeger constants are then defined as in the absolute case.

\mysubs{Cohomology expansion}
Let $\phi \in C^k(X)$ be a~$k$-cochain of $X$.

\begin{df}
The {\bf norm} of $\phi$ is
\[
\|\phi\|:=|\supp(\phi)|=|\{\sigma \in X(k): \ev{\phi}{\sigma} \neq 0\}|.
\]
The {\bf cosystolic norm} of $\phi$ is
\[
\|\phi\|_\csy:=\min_{\psi \in C^{k-1}(X)} \|\phi+d_{k-1}\psi\|.
\]
A cochain $\phi$ is a~{\bf cosystole} if $\|\phi\|=\|\phi\|_\csy$.
A~{\bf cosystolic form} of $\phi$ is any $\tilde \phi=\phi+d_{k-1}\psi$, such
that $\|\ti \phi\|=\|\phi\|_\csy$.
\end{df}

\begin{df}
The~{\bf coboundary expansion} of a~$k$-cochain $\phi \in C^k(X) \setminus B^k(X)$ is
\[\|\phi\|_{\exp}:=\|d_k \phi\|/\|\phi\|_\csy.\]
The {\bf $k$-th Cheeger constant} of $X$ is
\[
h^k(X):=\min_{\phi \in C^k(X) \setminus B_k(X)} \|\phi\|_{\exp}.
\]
\end{df}

Let us now turn to the relative case.
The $k$-th coboundary expansion of $k$-cochain $\phi \in C^k(X,Y)\sm
B^k(X,Y)$ is again
$\|\phi\|_{\exp}:=\|d_k \phi\|/\|\phi\|_{csy}$, and
the {\bf $k$-th Cheeger constant} of the pair $(X,Y)$ is:
\[
h^k(X,Y):=\min\{\|\phi\|_{\exp}:\phi \in C^k(X,Y)\sm B^k(X,Y)\}.
\]

Clearly, the non-relative norm, (co)systolic norm and (co)boundary
expansion are obtained by taking $Y$ to be the {\it void complex} $\emptyset=\{\,\}$  (see~\cite{book}), e.g., $h_k(X)=h_k(X,\emptyset)$,
$h^k(X)=h^k(X,\emptyset)$.

\begin{rem} $\,$

\noindent{\rm (1)} Let $c\in C_k(X,Y)$ and $\phi \in C^k(X,Y)$. Then
$\|c\|_{sys} \neq 0$ if and only if $c\not\in B_k(X,Y)$, and $\|\phi\|_{csy} \neq
0$ if and only if $\phi \not\in B^k(X,Y)$.

\noindent{\rm (2)} Note that $h_k(X,Y)>0$ if and only if $H_k(X,Y)=0$ and
$h^k(X,Y)>0$ if and only if $H^k(X,Y)=0$.  One can therefore view the expansion
constants $h_k(X,Y)$ and $h^k(X,Y)$ as refining the notion of
acyclicity, trying to catch phenomena which the regular (co)homology
does not. A possible analogy could be Whitehead torsion refining the
notion of homotopy equivalence.
\end{rem}

Here we study several aspects of the higher dimensional Cheeger constants.
The plan of the paper is as follows. In Section \ref{s:methods} we discuss some general tools that include:
\begin{itemize}
\item
A combinatorial lower bound on the cosystolic norm (Theorem \ref{thm:det1}).
\item
A lower bound on Cheeger constant using chain homotopy (Theorem \ref{t:lowbdh1}).
\item
An Alexander type duality between the homological and cohomological Cheeger constants (Theorem \ref{p:ad}).
\end{itemize}
Section \ref{s:sphere} is concerned with cosystoles and expansion of some concrete complexes and includes:
\begin{itemize}
\item A determination of the codimension one cosystoles and Cheeger constants of certain
$n$-pseudomanifolds (Theorem \ref{thm:sph}), and of Coxeter complexes (Theorem \ref{c:ecc}).
\item A lower bound on the Cheeger constants of a homogenous geometric lattices (Theorem \ref{t:explat}).
\end{itemize}
In Section \ref{s:cube} we revisit results of Gromov on expansion of products, including
\begin{itemize}
\item
An elementary derivation of Gromov's computation of the expansion of the hypercube (Theorem \ref{thm:cube}).
\item
A detailed proof of a theorem of Gromov on the expansion of a product of a complex with a simplex (Theorem \ref{t:prod}).
\end{itemize}

In Section \ref{s:maxs} we consider some extremal problems on cosystoles and expansion. These include
\begin{itemize}
\item
A lower bound on the maximal cosystole in a complex (Theorem \ref{p:lkn}).
\item
An upper bound on the expansion of bounded degree complexes (Theorem \ref{t:uphk}).
\item
A nearly sharp estimate on the cosystolic norm of the Paley cochains (Theorem \ref{p:paley}).
\end{itemize}
In Section \ref{s:nonab} we discuss
\begin{itemize}
\item
The non-abelian $1$-dimensional expansion of a simplex (Proposition \ref{bw1}).
\item
An application to a problem of Babson on bounded quotients of the fundamental group of a random $2$-complex (Theorem \ref{nonab}).
\end{itemize}
We conclude in Section \ref{s:conc} with some comments and open problems.

\section{General Tools}
\label{s:methods}

\mysubs{Detecting large cosystolic norm using cycles}
\label{subs:dct}

A question that frequently arises in specific examples, as well as in
theoretical context, is that of determining whether given cochain is a~cosystole.
Using the definition directly is impractical at best, since it would involve
going through all possible coboundaries and trying to see whether adding
one would reduce the norm of the cochain at hand.

The new idea which we introduce here is to use cycles to detect
in an indirect way that that our cocycle has a large cosystolic norm.
For this, we recall that coboundaries evaluate trivially on cycles, see (\ref{e:eval});
therefore, the evaluation of a cochain on a cycle does not change if we
add a coboundary to that cochain. In particular, if a cochain evaluates nontrivially
on a cycle, its support must intersect the support of that cycle, and that
will not change if we add a coboundary. The intersection cells may vary, but
the fact that the intersection is non-trivial will remain.

In its most basic form, our method is based on the fact that
if we have
a family of $t$ cycles with pairwise disjoint supports, and a cochain $\phi$ which evaluates
nontrivially on each of these cycles, then the cosystolic norm of $\phi$ is at least~$t$.
Let us now formalize these observations.
\begin{df}
Let $\cf \subset 2^V$ be a family of finite sets.
A subset $S \subset V$ is a {\bf piercing set} of $\cf$ if $S \cap F \neq \emptyset$ for all $F \in \cf$. The minimal cardinality of a piercing set of $\cf$, denoted by $\tau(\cf)$, is called the {\bf piercing number} of $\cf$.
\end{df}

\begin{thm}[The cycle detection theorem]
\label{thm:det1}
Let $X$ be a~polyhedral complex, and let $\phi$ be a~$k$-cochain of $X$.
Let $A=\{\alpha_1,\ldots,\alpha_t\}$ be a family of $k$-cycles of $X$, such that
$\ev{\phi}{\alpha_i} \neq 0$ for all $1 \leq i \leq t$. Let
$\cf=\{\supp(\alpha_1),\ldots,\supp(\alpha_t)\} \subset 2^{X(k)}$.
\begin{equation}
\label{e:cdt}
\|\phi\|_\csy \geq \tau(\cf).
\end{equation}
\end{thm}
\pr Let $\psi \in C^{k-1}(X)$. Then for any $1 \leq i \leq t$
\begin{equation*}\label{eq:det1}
\begin{split}
\ev{\phi+d_{k-1}\psi}{\alpha_i}&=\ev{\phi}{\alpha_i}+\ev{d_{k-1}\psi}{\alpha_i}\\
&=\ev{\phi}{\alpha_i}+\ev{\psi}{\bou{k}\alpha_i}\\
&=\ev{\phi}{\alpha_i}+\ev{\psi}{0}\\
&=\ev{\phi}{\alpha_i}\neq 0.
\end{split}
\end{equation*}
In particular, $\supp(\phi+d_{k-1}\psi)\cap\supp\alpha_i\neq\es$.
It follows that $\supp(\phi+d_{k-1}\psi)$ is a piercing set of $\cf$ and therefore
$\|\phi+d_{k-1}\psi\|=|\supp(\phi+d_{k-1}\psi)| \geq \tau(\cf)$. Since this is true for all $\psi$,
we get (\ref{e:cdt}).
\qed
\begin{crl} \label{thm:det2}
Let $X$ be a~polyhedral complex, and let $\phi$ be a~$k$-cochain of $X$.
Let $A=\{\alpha_1,\ldots,\alpha_t\}$ be a family of $k$-cycles of $X$ with pairwise disjoint supports, such that $\ev{\phi}{\alpha_i} \neq 0$ for all $1 \leq i \leq t$.
Then $\|\phi\|_{\csy} \geq t$.
\end{crl}
\noindent
{\bf Example:} Let $n=(k+2)m$ and let
$\dn$ denote the $(n-1)$-simplex on the vertex set $V=V_0 \cup \cdots \cup V_{k+1}$
where the $V_i$'s are disjoint of cardinality $m$.
Consider the collection of $k$-simplices
$S=\{[v_0,\ldots,v_k]: (v_0,\ldots,v_k) \in V_0 \times \cdots \times V_k\}$
and let $\phi=S^* \in C^k(\dn)$.
The following fact was mentioned in \cite{MW}.
\begin{claim}
\label{c:nk2}
\[
\|\phi\|_{csy}=\|\phi\|=m^{k+1}.
\]
\end{claim}
{\bf Proof:} For a $(k+2)$-tuple $\underline{v}=(v_0,\ldots,v_{k+1}) \in V_0 \times \cdots \times V_{k+1}$ let $\alpha_{\underline{v}}= \partial_{k+1}[v_0,\ldots,v_{k+1}] \in Z_k(\dn)$.
Identify each $V_i$ with a copy of the cyclic group $\Int_m$ and let
\[
T=\{(v_0,\ldots,v_{k+1}) \in V_0 \times \cdots \times V_{k+1}: v_0+\cdots+v_{k+1}=0\}.
\]
Let $A=\{\alpha_{\underline{v}}:\underline{v} \in T\}$. Clearly
$\supp(\alpha_{\underline{v}}) \cap \supp(\alpha_{\underline{u}})=\emptyset$ for
$\underline{u} \neq \underline{v} \in T$. Furthermore
$$\ev{\phi}{\alpha_{\underline{v}}}=\ev{\phi}{\partial_{k+1}[v_0,\ldots,v_{k+1}]}=
\ev{d_k\phi}{[v_0,\ldots,v_{k+1}]}=1.$$
Corollary \ref{thm:det2} therefore implies that
$\|\phi\|_{\csy} \geq |A|=m^{k+1}$.
\qed

\mysubs{Lower bounds for expansion via cochain homotopy}
\label{subs:lowbd1}

Let $X$ be an $n$-dimensional simplicial complex and let $k \leq n-1$.
Let $(S,\mu)$ be a finite probability space. Let $\left\{c_{s,\sigma}:(s,\sigma) \in S \times X(k)\right\}$
be a family of $(k+1)$-chains of $X$ and let
$\left\{c_{s,\tau}: (s,\tau) \in S \times X(k-1)\right\}$ be a family of $k$-chains of $X$,
such that for all $(s,\sigma) \in S \times X(k)$ we have
\begin{equation}
\label{e:fill}
\partial_{k+1} c_{s,\sigma}=\sigma+\sum_{j=0}^k c_{s,\sigma_j},
\end{equation}
where $\sigma_j$ denotes the $j$-th face of $\sigma$.
For $i=k,k+1$ and $s \in S$, define $T_s:C^i(X) \rightarrow C^{i-1}(X)$
as follows. For $\phi \in C^i(X)$ and $\sigma \in X(i-1)$ let
\begin{equation}
\label{e:defts}
\ev{T_s \phi}{\sigma}=\ev{\phi}{c_{s,\sigma}}.
\end{equation}

A natural way to think about \eqref{e:defts} and \eqref{e:fill} is to say
that the map $\sigma\mapsto c_{s,\sigma}$ is a~chain homotopy between the trivial
map and the identity map of $C_k(X)$, and that $T_s$ is its dual. It follows that $T_s$ satisfies
\begin{equation}
\label{e:contr}
d_{k-1}T_s+T_s d_k=\text{Id}_{C_k(X)}.
\end{equation}
For $(s,\sigma) \in S \times X(k)$ we set
\begin{equation}
\label{e:fsig}
F_{s,\sigma}:=\supp(c_{s,\sigma}) \subset X(k+1).
\end{equation}
For $\tau \in X(k+1)$ and $s \in S$ let
$\delta_s(\tau):=|\{\sigma \in X(k): \tau \in F_{s,\sigma}\}|$.
Let $E[f]$ denote the expectation of a random variable $f$ on $S$.
The following result is an elaboration of an idea from \cite{LMM}.
\begin{theorem}
\label{t:lowbdh1}
\begin{equation*}
\label{e:lowbdh1}
h^k(X) \geq \left( \max_{\tau \in X(k+1)} E\left[\delta_s(\tau)\right]\right)^{-1}.
\end{equation*}
\end{theorem}
\pr
Let $\phi \in C^k(X) \setminus B^k(X)$ and $s \in S$. By (\ref{e:contr})
\[
T_s d_k \phi=\phi -d_{k-1}T_s \phi,
\]
and hence
$\|\phi\|_{csy} \leq \|T_s d_k \phi\|$.
Taking expectation over $S$ we obtain
\begin{equation*}
\label{e:average}
\begin{split}
\|\phi\|_{csy}& \leq E\left[\|T_{s} d_k \phi\|\right]=\sum_{s \in S} \mu(s)\, \|T_{s} d_k \phi\| \\
&= \sum_{s \in S}\mu(s)\, |\{\sigma \in X(k): \ev{d_k\phi}{c_{s,\sigma}} \neq 0\}| \\
&= \sum_{s \in S}\mu(s)\, |\{\sigma\in X(k):|\supp(d_k \phi)\cap\supp(c_{s,\sigma})|
\textrm{ is odd}\}| \\
&\leq\sum_{s \in S} \mu(s)\,
|\{\sigma \in X(k): \supp(d_k \phi) \cap \supp(c_{s,\sigma}) \neq \emptyset\}| \\
&\leq\sum_{s\in S}\mu(s) \sum_{\sigma\in X(k)}|\supp(d_k \phi)\cap\supp(c_{s,\sigma})| \\
&= \sum_{s \in S} \mu(s)\sum_{\tau \in \supp(d_k\phi)}
|\{\sigma \in X(k): \tau \in \supp(c_{s,\sigma})\}| \\
&= \sum_{\tau \in \supp(d_k\phi)} \sum_{s \in S} \mu(s)\delta_s(\tau)= \sum_{\tau \in \supp(d_k\phi)} E\left[\delta_s(\tau)\right] \\
&\leq \|d_k \phi\| \cdot \max_{\tau \in X(k+1)} E\left[\delta_s(\tau)\right].
\end{split}
\end{equation*}
{\enp}
For $\tau \in X(k+1)$ let
\[
\delta(\tau):=\sum_{s \in S} \delta_s(\tau)=|\{(s,\sigma) \in S \times X(k): \tau \in \supp(c_{s,\sigma})\}|.
\]
Specializing Theorem \ref{t:lowbdh1} to the case of uniform distribution $\mu(s)\equiv\frac{1}{|S|}$,
we obtain the following
\begin{crl}
\label{c:lowbd}
\begin{equation*}
\label{e:lowbdh}
h^k(X) \geq \frac{|S|}{\max_{\tau \in X(k+1)} \delta(\tau)}.
\end{equation*}
\end{crl}

\mysubs{Alexander duality and expansion}
\label{ss:duality}
Let $\dn$ denote the $(n-1)$-simplex on an $n$-element vertex set $V$
and let $X$ be a~simplicial subcomplex of $\dn$. For a subset $\sigma
\subset V$ let $\overline{\sigma}=V-\sigma$. The Alexander dual $\xv$
of $X$ is the simplicial complex given by
\[\xv=\{\sigma \in \dn: \overline{\sigma}  \not\in X\}.\]
Note, that $(\xv)^\vee=X$. We also have $
\{\emptyset\}^\vee=\partial \dn \simeq S^{n-2}$, ${\rm
  void}^\vee=\dn$, $(\dn)^\vee={\rm void}$.

Let $Y$ be a subcomplex of $X$.  The combinatorial version of the
relative Alexander duality is the following
\begin{thm}[Alexander Duality]
\label{t:ad}
For $0 \leq k \leq n-1$
\[H_k(X,Y;\Int)\cong H^{n-k-2}(\yv,\xv;\Int).\]
\end{thm}
\noindent
In fact, there is a chain complex isomorphism $C_*(X,Y;G) \cong
C^*(Y^\vee,X^\vee;G)$, for an arbitrary abelian group~$G$.
The counterpart of Alexander duality for expansion is the following
\begin{theorem}
\label{p:ad}
For $0 \leq k \leq n-1$
\[h_k(X,Y)=h^{n-k-2}(\yv,\xv).\]
\end{theorem}

Letting $Y$ be the void simplicial complex in
Proposition~\ref{p:ad}, we obtain the following corollary.

\begin{crl}
Let $X \subset \dn$. Then:
\[h_k(X)=h^{n-k-2}(\dn,\xv),\]
\[h^k(X)=h_{n-k-2}(\dn,\xv).\]
\end{crl}

\prn{Theorem~\ref{p:ad}}
Define a~linear map $A_k:C_k(X,Y) \rightarrow C^{n-k-2}(\yv,\xv)$ as follows.
For a~generator $\sigma \in X(k)$ of $C_k(X,Y)$ and $\tau \in \yv(n-k-2)$ let
\[
\ev{A_k\sigma}{\tau}=\delta(\sigma,\overline{\tau})=\left\{
\begin{array}{ll}
1 & \tau=\overline{\sigma}, \\
0 & otherwise.
\end{array}
\right.~~\]
Note that $A_k$ is well-defined: If $\sigma \in X(k)$ and $\tau \in \xv(n-k-2)$
then $\ev{A_k\sigma}{\tau}=0$, thus $A_k\sigma \in C^{n-k-2}(\yv,\xv)$. Moreover, if
 $\sigma \in Y(k)$ then $\ev{A_k\sigma}{\tau}=0$ for all $\tau \in \yv(n-k-2)$, i.e.,
$A_k \sigma=0$. It is straightforward to check that $A_k$ is an isomorphism and
that it commutes with the differentials, i.e.,
\begin{equation}
\label{e:comm}
d_{n-k-2}A_k=A_{k-1}\partial_k.
\end{equation}
Observe that if $c=\sum_{\sigma \in X(k)} a_{\sigma} \sigma \in
C_k(X,Y)$ and $\tau \in \yv(n-k-2)$, then
$\ev{A_k c}{\tau}=a_{\overline{\tau}}$. Therefore
\begin{equation}
\label{e:du1}
\|A_k c\|=\|c\|.
\end{equation}
Combining (\ref{e:comm}) and (\ref{e:du1}) it follows that
\begin{equation}
\label{e:du2}
\begin{split}
\|A_k c\|_{csy}&=\min\{\|A_kc+d_{n-k-3}\psi\|:\psi \in C^{n-k-3}(\yv,\xv)\} \\
&=\min\{\|A_kc+d_{n-k-3}A_{k+1}c'\|:c' \in C_{k+1}(X,Y)\} \\
&=\min\{\|A_k c+A_k \partial_{k+1}c'\|: c' \in C_{k+1}(X,Y)\} \\
&=\min\{\|A_k(c+\partial_{k+1}c')\|: c' \in C_{k+1}(X,Y)\} \\
&=\min\{\|c+\partial_{k+1}c'\|: c' \in C_{k+1}(X,Y)\} \\
&=\|c\|_{sys}.
\end{split}
\end{equation}
Next note that (\ref{e:comm}) implies that $A_k$ maps
$C_k(X,Y)\setminus B_k(X,Y)$ injectively onto $C^{n-k-2}(\yv,\xv)
\setminus B^{n-k-2}(\yv,\xv)$. Therefore, if $c \in C_k(X,Y)
\setminus B_k(X,Y)$ then by (\ref{e:du1}) and (\ref{e:du2}):
\begin{equation}
\label{e:du3}
\begin{split}
h_k(c)&=\frac{\|\partial_k c\|}{\|c\|_{sys}} \\
&=\frac{\|A_{k-1} \partial_k c\|}{\|A_k c\|_{csy}} \\
&=\frac{\|d_{n-k-2} A_k c\|}{\|A_k c\|_{csy}} \\
&=h^{n-k-2}(A_k c).
\end{split}
\end{equation}
Theorem \ref{p:ad} now follows by minimizing (\ref{e:du3}) over all $c
\in C_k(X,Y) \setminus B_k(X,Y)$.  {\enp}


\section{Cosystoles and Expansion of Pseudomanifolds and Geometric Lattices}
\label{s:sphere}

\mysubs{The $(n-1)$-th Cheeger constant of an~$n$-pseudomanifold}

Let $X$ be an $n$-dimensional simplicial complex. The \emph{flip graph} of $X$, is the graph $G_X=(V_X,E_X)$ whose vertex set is $V_X=X(n)$ - the set of all
$n$-simplices of $X$, and whose edge set $E_X$ consists of all pairs
$\{\sigma,\sigma'\}$ such that $\dim (\sigma \cap \sigma')=n-1$.

Suppose now that $X$ is a triangulation of an $n$-pseudomanifold,
i.e., $X$ is a finite pure $n$-dimensional simplicial complex such
that any $\tau \in X(n-1)$ is contained in exactly two $n$-simplices
of~$X$ and such that $G_X$ is connected.  For $\phi \in C^{n-1}(X)$ let
$G_{\phi}=(V_X,E_{\phi})$ be the subgraph of $G_X$ with edge set
\[E_{\phi}=\big\{\{\sigma_1,\sigma_2\} \in E_X: \sigma_1 \cap \sigma_2
\in \supp \phi \big\}.\] Let $\cf_X$ denote the family of all
subgraphs $F=(V_X,E(F))$ of $G_X$ such that
\[|E(C) \cap E(F)| \leq \frac{|E(C)|}{2}\]
for any Eulerian subgraph $C=(V(C),E(C))$ of $G_X$. We note the
following properties of the family $\cf_X$.

\begin{claim}
\label{c:fx}
Let $F=(V_X,E(F)) \in \cf_X$. Then the following hold.
\begin{enumerate}
\item[(i)] The graph $F$ is a forest.
\item[(ii)] If vertices $u,v \in V_X$ are in the same tree component
  of $F$ then we have $\dist_F(u,v)=\dist_{G_X}(u,v)$. In particular,
  if $P=(V(P),E(P))$ is a~path in~$F$, then $|E(P)| \leq \diam(G_X)$.
\end{enumerate}
\end{claim}

\pr To prove (i) note that if $C=(V(C),E(C))$ is a cycle in $G_X$ then
\[|E(C) \cap E(F)| \leq \frac{|E(C)|}{2}<|E(C)|,\] so in particular
$E(C) \not\subset E(F)$. To show (ii), let $P=(V(P),E(P))$ be the path
in $F$ between $u$ and $v$ and let $Q=(V(Q),E(Q))$ be a minimal $u-v$
path in $G_X$. Consider the Eulerian graph $R=(V_X,E(R))$ where
$E(R)=(E(P)\setminus E(Q)) \cup (E(Q)\setminus E(P))$. Then
\begin{equation}
\label{e:pqr}
\begin{split}
|E(P)|&=|E(P) \cap E(Q)|+|E(P)\setminus E(Q)| \\
&\leq |E(P) \cap E(Q)|+|E(F) \cap E(R)| \\
&\leq |E(P) \cap E(Q)|+\frac{|E(R)|}{2} \\
&=|E(P) \cap E(Q)|+\frac{|E(P)\setminus E(Q)|}{2}+\frac{|E(Q)\setminus E(P)|}{2} \\
&=\frac{1}{2} (|E(P)|+|E(Q)|).
\end{split}
\end{equation}
Hence $\dist_F(u,v)=|E(P)| \leq |E(Q)|=\dist_{G_X}(u,v)$.
\qed

Using Claim \ref{c:fx} we next give a combinatorial description of the
$(n-1)$-cosystoles in $X$.

\begin{claim}
\label{c:dualg}
Let $X$ be an~arbitrary $n$-pseudomanifold and let $\phi \in C^{n-1}(X)$.
Then the following hold.
\begin{enumerate}
\item[(i)] The mapping $\phi \rightarrow G_{\phi}$ maps $C^{n-1}(X)-\{0\}$
  injectively onto all subgraphs of the graph~$G_X$.
\item[(ii)] We have
\begin{equation}
\label{e:expodd}
\supp(\cob{n-1}\phi)=\{\sigma \in V_{\phi}: \deg_{G_{\phi}}(\sigma) ~odd~\}.
\end{equation}
\item[(iii)] Suppose $H^{n-1}(X;\Int_2)=0$. Then $\phi$ is a cosystole
 if and only if $G_{\phi} \in \cf_X$.
\end{enumerate}
\end{claim}
\pr Parts (i) and (ii) follow directly from the definitions. We
proceed to prove part~(iii).

First suppose that $\|\phi\|=\|\phi\|_{csy}$. Assume $C=(V_X,E(C))$ is an
Eulerian subgraph of $G_{\phi}$ with edge set
\[
E(C)=\big\{\{\sigma_0,\sigma_1\},\{\sigma_1,\sigma_2\},\ldots,\{\sigma_{m-1},\sigma_m\}\big\}.
\]
Set
$\psi:=\sum_{i=1}^m(\sigma_{i-1}\cap\sigma_i)^*$.  Clearly $\psi$ is
a~cocycle. On the other hand, we assumed that $H^{n-1}(X;\Int_2)=0$,
so $\psi$ must also be a~coboundary. We there\-fore have
\begin{equation*}\label{e:ephi}
\begin{split}
2|E_{\phi} \cap E(C)|&=2|E_{\phi} \cap E_{\psi}| =|E_{\psi}|+|E_{\phi}|-|E_{\phi+\psi}| \\
&=\|\psi\|+\big(\|\phi\|-\|\phi+\psi\|\big)\leq \|\psi\|=|E(C)|.
\end{split}
\end{equation*}
Conversely, suppose that $G_{\phi} \in \cf_X$ and let $\psi\in
B^{n-1}(X)$. Then $G_{\psi}$ is Eulerian and hence $|E_{\phi} \cap E_{\psi}|
\leq{|E_{\psi}|}/{2}$. Therefore
\begin{equation}
\label{e:euler}
\begin{split}
\|\phi+\psi\|&=|E_{\phi+\psi}|=|E_{\phi}|+|E_{\psi}|-2|E_{\phi} \cap E_{\psi}| \\
&\geq |E_{\phi}|=\|\phi\|.
\end{split}
\end{equation}
We conclude that $\|\phi\|_{csy}=\|\phi\|$.
{\enp}
\noindent
Claim \ref{c:dualg} implies the following combinatorial
characterization of the $(n-1)$-coboundary expansion of
$n$-pseudomanifolds. See Lemmas 2.4 and 2.5 in \cite{SKM} for a
related result.
\begin{thm}
\label{thm:sph}
Let $X$ be an $n$-pseudomanifold such that $H^{n-1}(X;\Int_2)=0$. Then
\begin{equation}
\label{eq:pseh}
h^{n-1}(X)=\frac{2}{\diam(G_X)}.
\end{equation}
\end{thm}
\pr
In view of Claim \ref{c:dualg}, it suffices to show that
\begin{equation}
\label{e:graphv}
\min_{F=(V,E)\in \cf_X}\frac{\big|\{v: \deg_F(v)\text{ is odd }\}\big|}
{|E|}=\frac{2}{\diam(G_X)}.
\end{equation}
To prove the lower bound let $F=(V_X,E(F)) \in \cf_X$ and set
\[k:=\big|\{v \in V_X: \deg_F(v) \text{ is odd }\}\big|/2.\]
Claim \ref{c:fx}(i) implies that $F$ is a~forest. It follows (see,
e.g., Theorem 2.1.10 in \cite{West}) that there exist $k$ edge disjoint
paths $P_1=(V_1,E_1),\ldots,P_k=(V_k,E_k)$ in $F$ such that $E_1\cup
\cdots \cup E_k=E(F)$.  Claim \ref{c:fx}(ii) implies that $|E_i| \leq
\diam(G_X)$, hence
\begin{equation}
\label{e:lowbd}
\begin{split}
\frac{\big|\{v: \deg_F(v) ~\text{odd}\}\big|}{|E(F)|}&=\frac{2k}{\sum_{i=1}^{k} |E_i|} \\
&\geq \frac{2k}{k \cdot \diam(G_X)} =\frac{2}{\diam(G_X)}.
\end{split}
\end{equation}
Finally, we show that equality in (\ref{e:graphv}) is attained for
some $F \in \cf_X$. Let $u,v \in V_X$ such that
$\dist_{G_X}(u,v)=\diam(G_X)$ and let $P=(V_X,E(P))$ be a minimal
$u-v$ path in $G_X$. Clearly $P \in \cf_X$ and
\[
\frac{\big|\{v: \deg_P(v) ~\text{odd}\}\big|}{|E(P)|}=\frac{2}{\diam(G_X)}.
\]
This shows \eqref{eq:pseh}.
{\enp}

\mysubs{The Expansion of Coxeter Complexes}
\label{subs:cox}
Let $W$ be an arbitrary finite Coxeter group with the set of
generators $S$ and a root system~$\Phi$. We refer to the books by
Humphrey \cite{Hu} and by Ronan \cite{Ro} for the
theory of Coxeter groups and Coxeter complexes. For $J\subset S$ let
$W_J=\langle s:s\in J\rangle$ be the subgroup of $W$ generated by $J$.
For $s\in S$ let $(s)=S-\{s\}$.

\begin{df}
The \emph{Coxeter complex} $\Delta(W,S)$ is the simplicial complex on the vertex set
$V=\bigcup_{s \in S} W/W_{(s)}$ whose maximal simplices are
$C_w=\{wW_{(s)}:s \in S\}$, for $w \in W$.
\end{df}

The simplicial complex $\Delta(W,S)$ is
a~triangulation of $(|S|-1)$-dimensional sphere. It is well-known, see, e.g.,
\cite[Theorem 2.15]{Ro}, that $\diam(G_{\Delta(W,S)})={|\Phi|}/{2}$.
Therefore by Theorem~\ref{thm:sph} we have
\begin{crl}
\label{c:ecc}
$h^{|S|-2}\left(\Delta(W,S)\right)={4}/{|\Phi|}$.
\end{crl}
\noindent
{\bf Examples:} \\ (i) Let $W=\ccs_n$ be the symmetric group on
$[n]$ with the set of generators $S=\{s_1,\ldots,s_{n-1}\}$ where
$s_i=(i,i+1)$ for $1 \leq i \leq n-1$.  Then $|\Phi|=n(n-1)$ and
$\Delta(W,S)$ is isomorphic to $\text{sd} \, \partial\smp{n-1}$, the
barycentric subdivision of the boundary of the $(n-1)$-simplex. Hence
\begin{equation}
\label{e:sdsim}
h^{n-3}\left(\text{sd} \ \partial\smp{n-1}\right)=\frac{4}{n(n-1)}.
\end{equation}
We next describe an explicit $(n-3)$-cochain $\phi_n$ of $X_n:=\text{sd} \, \partial\smp{n-1}$
such that $\|\phi_n\|_{\exp}=\frac{4}{n(n-1)}$.
With a permutation $\pi=(\pi(1),\ldots,\pi(n)) \in \ccs_n$ we associate the $(n-2)$-face $F(\pi)$ of $X_n$ given by
$$F(\pi)=\big[\{\pi(1)\} \subset \{\pi(1),\pi(2)\} \subset \cdots \subset \{\pi(1),\ldots,\pi(n-1)\}\big].
$$
For $1 \leq i \leq n-1$, the $i$-th face of $F(\pi)$ is:
$$F(\pi)_i= F(\pi) \setminus \big\{\,\{\pi(1),\ldots,\pi(i)\}\,\big\}.$$
Define a sequence of permutations $\pi_0,\ldots,\pi_{\binom{n}{2}} \in \ccs_n$ as follows.
First let $\pi_0=(1,\cdots,n)$ be the identity permutation.
Next let $1 \leq m \leq \binom{n}{2}$. Then $m$ can be written uniquely
as
$$m=m(j,\ell):=(j-1)n-\binom{j}{2}+\ell$$
where $1 \leq j \leq n-1$ and $1 \leq \ell \leq n-j$.
Define
$$
\pi_{m(j,\ell)}=(n,n-1,\ldots,n-j+2,1,2,\ldots,n-j-\ell,n-j+1,n-j-\ell+1,\cdots,n-j).
$$
Let
$$\phi_n =\sum_{j=1}^{n-1}\sum_{\ell=1}^{n-j} F(\pi_{m(j,\ell)})_{n-\ell}^*
\in C^{n-3}(X_n).$$
It is straightforward to check that $F(\pi_m)_{n-\ell}=F(\pi_{m-1})_{n-\ell}$.
Hence $$\supp\big( d_{n-3} F(\pi_{m})_{n-\ell}^*\big)
=\{F(\pi_m),F(\pi_{m-1})\}$$
and therefore
$$\supp\left(d_{n-3}(\phi_n)\right)=\{F(\pi_0),F(\pi_{\binom{n}{2}})\}.$$
It can also be shown that $\phi_n$ is an $(n-3)$-cosystole, i.e., $\|\phi_n\|_{\csy}=\|\phi_n\|=\binom{n}{2}$.
It thus follows that
$$
\|\phi_n\|_{\exp}=\frac{\|d_{n-3} \phi_n\|}{\|\phi_n\|_{\csy}}=\frac{4}{n(n-1)}.
$$

\noindent
(ii) Let $W=\ccs_2 \wr \ccs_n$ be the hyperoctahedral group with the set
of generators $S=\{\epsilon,s_1,\ldots,s_{n-1}\}$ where $\epsilon=(1,2)
\in \ccs_2$ and $s_i=(i,i+1) \in \ccs_n$ for $1 \leq i \leq n-1$.  Then
$|\Phi|=2n^2$ and $\Delta(W,S)$ is isomorphic to $\text{sd} \,
(\partial\smp{1})^{*n}$, the barycentric subdivision of the
octahedral $(n-1)$-sphere.  Hence
\begin{equation}
\label{e:sdocta}
h^{n-2}\left( \text{sd} \, (\partial\smp{1})^{*n} \right)=\frac{2}{n^2}.
\end{equation}

\editcut

\mysubs{Expansion of Homogenous Geometric Lattices}
\label{subs:geolat}

Let $(P,\leq)$ be a finite poset. The \emph{order complex} of $P$ is
the simplicial complex on the vertex set $P$ whose simplices are the
chains $a_0<\cdots < a_k$ of $P$, see \cite{book}. In the sequel we identify a poset
with its order complex.
A poset $(L,\leq)$ is a \emph{lattice} if any two elements $x,y \in L$ have a unique minimal upper bound $x \vee y$ and a unique maximal lower bound $x \wedge y$. A lattice $L$ with minimal element $\hzero$ and maximal element $\hone$ is \emph{ranked}, with rank function
$\rk(\cdot)$, if $\rk(\hzero)=0$ and $\rk(y)=\rk(x)+1$ whenever $y$ covers $x$. $L$ is a \emph{geometric lattice} if $\rk(x)+\rk(y) \geq \rk(x \vee y)+\rk(x \wedge y)$ for any $x,y \in L$,
and any element in $L$ is a join of atoms (i.e., rank $1$ elements).

Let $L$ be a geometric lattice of rank $\rk(\hone)=n$, and let $\overline{L}=L-\{\hzero,\hone\}$. A classical result of Folkman \cite{F66} asserts that $\thh_{i}(\overline{L})=0$ for $i < n-2$.
It is thus natural to ask for lower bounds on the Cheeger constants $h^i(\overline{L})$ for $i<n-2$.
In this section we approach this question using the cochain homotopy method of Section \ref{subs:lowbd1}. Let $S$ be a set of linear orderings on the set of atoms $A$ of $L$.
Let $\prec_{s}$ denote the ordering associated with $s \in S$.
For a subset $\{b_1,\ldots,b_m\} \subset A$ such that $m \leq n-1$ let
\[
K(b_1,\ldots,b_m)=\sum_{\pi \in \ccs_m} [b_{\pi(1)},b_{\pi(1)} \vee b_{\pi(2)},\ldots, b_{\pi(1)} \vee b_{\pi(2)} \vee \cdots \vee b_{\pi(m)}] \in C_{m-1}(\overline{L}).
\]
Note that
\begin{equation}
\label{e:bndb}
\partial_{m-1} K(b_1,\ldots,b_m)=\sum_{i=1}^m K(b_1,\ldots,\widehat{b_i},\ldots,b_m).
\end{equation}
Let $-1 \leq k \leq n-3$ and let $\sigma=[v_0<\cdots<v_k]$ be a $k$-simplex of $\overline{L}$.
Fix $s \in S$. Let $a_{s,k+1}(\sigma)=\min A$, and for $0 \leq i \leq k$ let
$a_{s,i}(\sigma)=\min\{a \in A: a \leq v_i\}$, where
both minima are taken with respect to $\prec_{s}$.
Note that
$$
a_{s,k+1}(\sigma) \preceq_{s} \cdots \preceq_{s} a_{s,0}(\sigma).
$$
Define
$$
c_{s,\sigma}=\sum_{j=0}^{k+1} K\left(a_{s,0}(\sigma),\ldots,a_{s,j}(\sigma)\right)*[v_j,\ldots,v_k],
$$
where $[v_j,\ldots,v_k]$ is interpreted as the empty simplex if $j=k+1$.
Note that
\begin{equation}
|\supp(c_{s,\sigma})| \leq \sum_{j=1}^{k+2} j!.
\end{equation}

\begin{claim}
\label{c:jchk}
For all $s \in S$, $0 \leq k \leq n-3$ and $\sigma \in \overline{L}(k)$
\begin{equation}
\label{e:chom}
\partial_{k+1} c_{s,\sigma}=\sigma+\sum_{i=0}^k c_{s,\sigma_i}.
\end{equation}
\end{claim}
\noindent
{\bf Proof:} As $s$ is fixed, we abbreviate $c_{\sigma}=c_{s,\sigma}$ and $a_j=a_{s,j}$.
Let $0 \leq i \leq k$ then
\begin{equation}
\label{e:bdcs}
\begin{split}
c_{\sigma_i}=c_{[v_0,\ldots,\widehat{v_i},\ldots,v_k]}
&=\sum_{j=0}^{i-1} K(a_0,\ldots,a_j)*[v_j,\ldots,\widehat{v_i},\ldots,v_k] \\
&+\sum_{j=i+1}^{k+1}K(a_0,\ldots,\widehat{a_i},\ldots,a_j)*[v_j,\ldots,v_k].
\end{split}
\end{equation}
Hence
\begin{equation}
\label{e:bdcs0}
\begin{split}
\sum_{i=0}^{k} c_{\sigma_i}
&=\sum_{i=0}^{k}\sum_{j=0}^{i-1} K(a_0,\ldots,a_j)*[v_j,\ldots,\widehat{v_i},\ldots,v_k] \\
&+\sum_{i=0}^{k}\sum_{j=i+1}^{k+1}K(a_0,\ldots,\widehat{a_i},\ldots,a_j)*[v_j,\ldots,v_k].
\end{split}
\end{equation}
Using (\ref{e:bndb}) and (\ref{e:bdcs0}) we compute
\begin{equation}
\label{e:bdcs1}
\begin{split}
\partial_{k+1}c_{\sigma} &= \sum_{j=0}^{k+1} \partial_{j} K(a_0,\ldots,a_j)*[v_j,\ldots,v_k] \\
&+ \sum_{j=0}^{k} K(a_0,\ldots,a_j)*\partial_{k-j}[v_j,\ldots,v_k] \\
&=\sum_{j=0}^{k+1}\sum_{i=0}^j
K(a_0,\ldots,\widehat{a_i},\ldots,a_j)*[v_j,\ldots,v_k] \\
&+ \sum_{j=0}^k \sum_{i=j}^k
K(a_0,\ldots,a_j)*[v_j,\ldots,\widehat{v_i},\ldots,v_k] \\
&=\sum_{i=0}^{k+1}\sum_{j=i}^{k+1}
K(a_0,\ldots,\widehat{a_i},\ldots,a_j)*[v_j,\ldots,v_k] \\
&+ \sum_{i=0}^k \sum_{j=0}^i
K(a_0,\ldots,a_j)*[v_j,\ldots,\widehat{v_i},\ldots,v_k] \\
&=\sum_{i=0}^{k}\sum_{j=i+1}^{k+1}
K(a_0,\ldots,\widehat{a_i},\ldots,a_j)*[v_j,\ldots,v_k]  \\
&\,\,\,\,\,\,\,\,\,+ \sum_{i=0}^k K(a_0,\ldots,a_i)*[v_{i+1},\ldots,v_k] +[v_0,\ldots,v_k] \\
&+ \sum_{i=0}^k \sum_{j=0}^{i-1}
K(a_0,\ldots,a_j)*[v_j,\ldots,\widehat{v_i},\ldots,v_k] \\
&\,\,\,\,\,\,\,\,\,+\sum_{i=0}^k K(a_0,\ldots,a_i)*[v_{i+1},\ldots,v_k] \\
&=\sigma+\sum_{i=0}^k c_{\sigma_i}.
\end{split}
\end{equation}
{\enp}

\begin{figure}
\begin{center}
  \scalebox{0.4}{\input{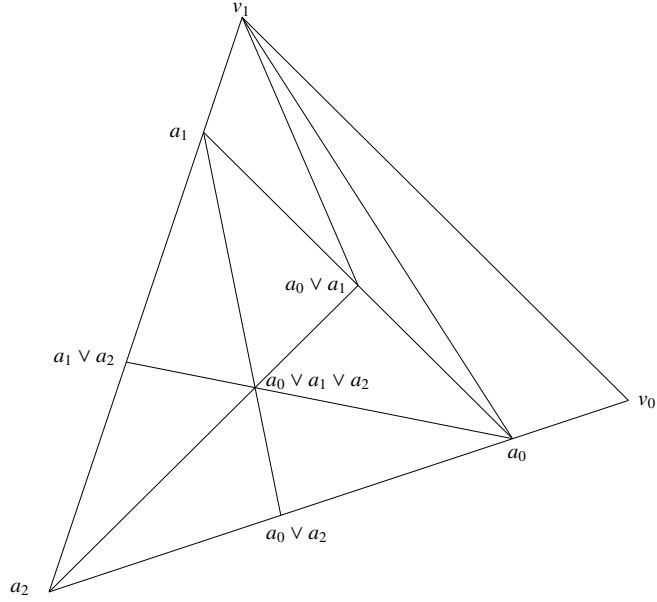}}
  \caption{$c_{s,\sigma}$ for $\sigma=[v_0,v_1]$.}
  \label{figure1}
\end{center}
\end{figure}

One natural choice of a set $S$ of linear orderings is the following.
Let $\prec$ be an arbitrary fixed linear order on the set of atoms $A$.
Let $S=Aut(L)$ be the automorphism group of $L$. For $s \in S$ let $\prec_s$ be the linear order
on $A$ defined by $a \prec_s a'$ if $s^{-1}(a) \prec s^{-1}(a')$.
Let $\idd$ denote the identity element of $S$. It is straightforward to check that
if $\sigma \in \overline{L}(k)$ and $0 \leq i \leq k+1$, then $a_{\idd,i}(s^{-1}(\sigma))=s^{-1}(a_{s,i}(\sigma))$ and hence
$c_{s,\sigma}=s\left(c_{\idd,s^{-1}(\sigma)}\right)$.
Using the definition of $F_{s,\sigma}$ (see Subsection \ref{subs:lowbd1} and in particular Eq. (\ref{e:fsig}) ), it follows that for any $s,t \in S$, $\sigma \in \overline{L}(k)$ and $\tau \in \overline{L}(k+1)$, it holds that
$\tau \in F_{s,\sigma}$ if and only if $t(\tau) \in F_{ts,t(\sigma)}$. In particular,
\begin{equation}
\label{e:degtau}
\delta(\tau)=\delta(t(\tau)).
\end{equation}

\begin{df}
\label{d:homl}
A geometric lattice $L$ is \emph{homogenous} if its automorphism group $G=Aut(L)$ is transitive on the set $\overline{L}(n-2)$ of top dimensional simplices of $\overline{L}$.
\end{df}
\begin{theorem}
\label{t:explat}
If $L$ is a homogenous geometric lattice of rank $n$ then
\begin{equation}
\label{e:explat}
h^{n-3}\left(\overline{L}\right) \geq \frac{f_{n-2}\left(\overline{L}\right)}{f_{n-3}\left(\overline{L}\right) \sum_{j=1}^{n-1}j!}.
\end{equation}
\end{theorem}
\noindent
{\bf Proof:} The homogeneity of $L$ together with (\ref{e:degtau}) imply that $\delta(\tau)=D$ is constant for all $\tau \in \overline{L}(n-2)$. Therefore
\begin{equation}
\label{e:vald}
\begin{split}
D\cdot f_{n-2}(\overline{L})&=\sum_{s \in S}\sum_{\sigma \in \overline{L}(n-3)} |F_{s,\sigma}| \\
&\leq |S| \cdot f_{n-3}(\overline{L})\cdot \sum_{j=1}^{n-1} j!.
\end{split}
\end{equation}
Hence, by Corollary \ref{c:lowbd}
\[
h^{n-3}\left(\overline{L}\right) \geq \frac{|S|}{D} \geq
\frac{f_{n-2}\left(\overline{L}\right)}{f_{n-3}\left(\overline{L}\right) \sum_{j=1}^{n-1}j!}.
\]
{\enp}
The spherical building $A_{n-1}(\FF_q)$ is the order complex $\overline{L}$,
where $L$ is the lattice of linear subspaces of $\FF_q^n$.
In \cite{Gr,LMM} it is shown that $h^{n-3}\left(A_{n-1}(\FF_q)\right) \geq \frac{q+1}{(n-1)n!}$.
Applying Theorem \ref{t:explat} to $A_{n-1}(\FF_q)$ and
noting that $f_{n-2}\left(A_{n-1}(\FF_q)\right) \cdot(n-1)=f_{n-3}\left(A_{n-1}(\FF_q)\right)\cdot (q+1)$, we obtain the following slight improvement.
\begin{crl}
\begin{equation*}
h^{n-3}\left(A_{n-1}(\FF_q)\right) \geq
\frac{f_{n-2}\left(A_{n-1}(\FF_q)\right)}{f_{n-3}\left(A_{n-1}(\FF_q)\right) \sum_{j=1}^{n-1}j!}=
\frac{q+1}{(n-1)\sum_{j=1}^{n-1}j!}.
\end{equation*}
\end{crl}
\editcut

\section{Coboundary expansion and products}
\label{s:cube}

\mysubs{Hypercube}
Let $Q_d$ denote the hypercube in dimension $d\geq 2$.  For future
reference note that number of vertices of $Q_d$ is $2^d$, number of
edges is $d\cdot 2^{d-1}$, etc; in general number of $k$-dimensional
cells is $\binom{d}{k}\cdot 2^{d-k}$.  The $k$-dimensional cells are
indexed by $d$-tuples of symbols $\{+,-,*\}$, where the total number
of occurences of $*$ is $k$.

It is certainly well-known that $h_0(Q_d)=1$.  Still, here is an
elementary argument.  Note that $h_0(Q_d)=1$ simply says that if $S$
is any set of vertices such that $|S|\leq 2^{d-1}$, then at least
$|S|$ edges connect $S$ to its complement. The equality is achieved if
for example $S$ consists of all vertices with the first coordinate
$+$. We can now easily prove this statement by induction on $d$. The
base $d=2$ is clear. For the induction step, let $V_+$ denote the set
of vertices of $Q_d$ with the first coordinate $+$ and let $V_-$
denote the set of vertices of $Q_d$ with the first coordinate
$-$. Accordingly set $S_+:=S\cap V_+$, $S_-:=S\cap V_-$, and
$T_+:=V_+\sm S_+$, $T_-:=V_-\sm S_-$.  Without loss of generality
assume that $|S_+|\leq|S_-|$. Let $e_{ij}$ denote the number of edges
between $S_i$ and $T_j$, for $i,j\in\{+,-\}$, and let $e$ denote the
number of edges between $S$ and its complement.  By induction
assumption, we have $e_{++}\geq |S_+|$. If also $|S_-|\leq |T_-|$, we
can apply induction assumption to $S_-$ as well; so we get
$e_{--}\geq|S_-|$ and are done.  Assume then we have $|S_-|\geq
|T_-|$, in which case we have $e_{--}\geq |T_-|$.  We have $e_{-+}\geq
|S_-|-|S_+|$, because each vertex in $S_-$ is connected by an edge to
exactly one vertex in $V_+$.  In total, we have $e\geq
e_{++}+e_{--}+e_{+-}\geq |S_+|+|T_-|+|S_-|-|S_+|=2^{d-1}$, and we are
done.

\begin{lm}
$h_k(Q_d)\leq 1$, for all $0\leq k\leq d-1$.
\end{lm}
\pr Fix $k\geq 1$, and let $E_k$ be the set of all $k$-cubes
indexed by $d$-tuples $(\bar x,*,\dots,*,0)$, where $\bar x$ is an
arbitrary $(d-k-1)$-tuple of $\{+,-\}$; so the number of $*$'s is $k$.
Clearly, $|E_k|=2^{d-k-1}$. Let $E^*_k$ be the $k$-cochain obtained by
summing up the characteristic cochains of the elements of $E_k$.

We can calculate the cosystole of $E_K^*$ by using our {\it detecting
  cycles} argument, see Section \ref{subs:dct}. As detecting cycles we take the boundaries of
$(\bar x,*,\dots,*)$, where $\bar x$ is an arbitrary $(d-k-1)$-tuple
of $\{+,-\}$. It will show that $\csy(E_k^*)=2^{d-k-1}$.  On the other
hand, an easy calculation shows that $|d_k(E_k^*)|=2^{d-k-1}$, so we
get $h_k(Q_d)\leq 1$. \enp

The lower bound can be shown by elementary methods as well.  Before we
proceed with the proof, let us show the following
inequality.

\begin{figure}[hbt]

  \input{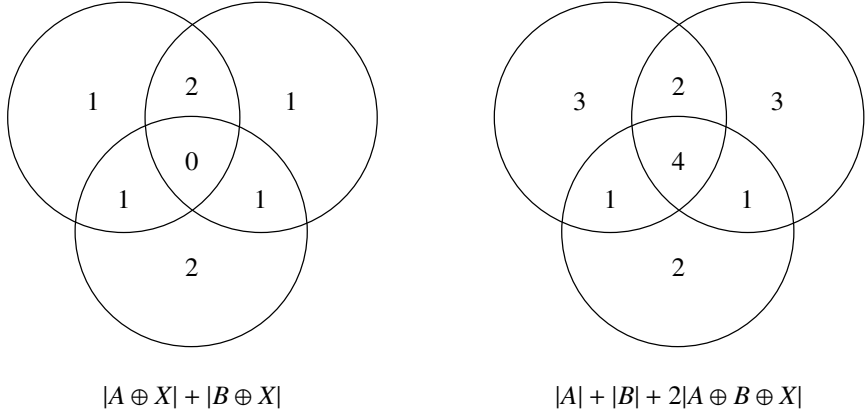}

\caption{Venn diagram illustrating the proof of Lemma~\ref{lm:axbx}.}
\label{fig:axbx}
\end{figure}

\begin{lm} \label{lm:axbx}
Let $A$, $B$, and $X$ be subsets of some universal set,
then we have
\begin{equation}
\label{eq:setin}
\sn{A\xor X}+\sn{B\xor X}\leq\sn{A}+\sn{B}+2\sn{A\xor B\xor X}.
\end{equation}
\end{lm}
\pr
Indeed, the inequality \eqref{eq:setin} follows from the the following calculation
\[
\begin{split}
\sn{A\xor X}+\sn{B\xor X}&=\sn{A}+\sn{X}-2\sn{A\cap X}+\sn{B}+\sn{X}-2\sn{B\cap X}\\
&=\sn{A}+\sn{B}+2\left(\sn{X}-\sn{A\cap X}-\sn{B\cap X}\right)\\
&\leq \sn{A}+\sn{B}+2\sn{X\sm(A\cup B)}\\
&\leq\sn{A}+\sn{B}+2\sn{A\xor B\xor X},
\end{split}\]
see Figure~\ref{fig:axbx}.
\qed

\begin{thm}\label{thm:cube}
We have $h_k(Q_d)=1$, for all $d\geq 2$, $0\leq k\leq d-1$.
\end{thm}
\pr We just need to show that $h_k(Q_d)\geq 1$. In other words, for
any cochain $\fy$, which is not a~cocycle, we need to show that
\begin{equation}\label{eq:ineq1}
\csyn{\fy}\leq\n{d^*\fy}.
\end{equation}
Our proof goes by induction on $d$ and $k$.  We already know
that~\eqref{eq:ineq1} holds for $k=0$ and arbitrary $d$. Furthermore,
when $k=d-1$, all the cochains which are not cocycles have both the
cosystolic norm as well as the norm of the coboundary equal to~$1$,
so~\eqref{eq:ineq1} becomes an~equality. This gives the boundary conditions
for the induction. To prove that \eqref{eq:ineq1} holds for $(d,k)$ we will
use the induction assumption that it holds for $(,)$.

Let us set $P_0:=\mycb{(a_1,\dots,a_t)\,|\,a_1=0}$, and
$P_1:=\mycb{(a_1,\dots,a_t)\,|\,a_1=1}$. Set furthermore
$P_*:=\mycb{(a_1,\dots,a_t)\,|\,a_1=*}$. This means that we fix one of
the directions of the hypercube and break the cube into two identical
lower-dimensional copies, called $P_0$ and $P_1$. The set $P_*$ contains
all the cubes which span across between the two halves. This in itself
is of course not a~subcomplex.

Let $\fy$ be an arbitrary cochain. Set $S_0:=\supp\fy\cap P_0$,
$S_1:=\supp\fy\cap P_1$, $S_*:=\supp\fy\cap P_*$, and
$\fy_0:=\fy_{S_0}$, $\fy_1:=\fy_{S_1}$, $\fy_*:=\fy_{S_*}$. Then $\fy$
decomposes as a~sum $\fy=\fy_0+\fy_1+\fy_*$, and furthermore, we have
\begin{equation}\label{eq:q11}
\n{\fy}=\n{\fy}_0+\n{\fy_1}+\n{\fy_*}.
\end{equation}
Assume $\alpha_0$ is a~cochain, such that
\begin{enumerate}
\item[(1)] $\supp\alpha_0\in P_0$;
\item[(2)] $\n{\alpha_0}=\csyn{\alpha_0}=\csyn{\fy_0}$, where
the cosystolic norm is taken
in~$P_0$;
\item[(3)] $\alpha_0-\fy_0=d^0\beta_0$, for some $\beta_0$, where
$d^0$ denotes coboundary operator in~$P_0$.
\end{enumerate}
Let us consider what happens if we replace $\fy$ with $\ti\fy=\fy+d\beta_0$,
where the coboundary operator is just the usual one in $X$. Since the
cochain is changed by a~coboundary, we have $\csyn{\fy}=\csyn{\ti\fy}$
and also $d(\ti\fy)=d(\fy)$. So if we show $\csyn{\ti\fy}\leq\n{d\ti\fy}$,
then we also show that $\csyn{\fy}\leq\n{d\fy}$. Set
$T:=\supp\ti\fy\cap P_0=\supp(\fy+d\beta_0)\cap P_0$. Clearly, we then have
\[T=\supp(\fy_0+d^0\beta_0)=\supp\alpha_0.\]
So $\ti\fy_T=\alpha_0$. Thus, replacing $\fy$ with $\ti\fy$ simply makes sure
that $\n{\fy_0}=\csyn{\fy_0}$, without changing either the cosystolic
norm or its coboundary. Completely identical argument holds for~$P_1$.

Summarizing our argument we conclude that we can add boundaries to
$\fy$ to make sure that $\n{\fy_0}=\csyn{\fy_0}$ and
$\n{\fy_1}=\csyn{\fy_1}$, where the cosystolic norm is taken in the
lower-dimensional cubes $P_0$ and $P_1$. By induction assumption for
$(d-1,k)$, we can therefore assume that
\begin{equation}
\label{eq:q0}
\n{\fy_0}\leq\n{d^0\fy_0}\textrm{ and }\n{\fy_1}\leq\n{d^1\fy_1},
\end{equation}
where $d^0$ and $d^1$ are the coboundary operators in $P_0$ and~$P_1$.

Let now $p_0\fy_*$ be the cochain in $P_0$ obtained from $\fy_*$ by
replacing the $*$ by $0$ in the first coordinate. Clearly, we have
\[d(p_0\fy_*)=\fy_*+d^0(p_0\fy_*).\]
Same way, let $p_1\fy_*$ be the cochain in $P_1$ obtained from $\fy_*$
by replacing the $*$ by~$1$ in the first coordinate. Again, we have
$d(p_1\fy_*)=\fy_*+d^1(p_1\fy_*)$. Finally, let $p_0(d\fy_*)$,
resp.\ $p_1(d\fy_*)$, denote the cochain in $P_0$, resp.\ $P_1$,
obtained from $d\fy_*$ by replacing the $*$ by $0$, resp.\ $1$, in the
first coordinate.  We have the following inequalities:
\begin{equation}
\begin{split}
\label{eq:q1}
2\cdot\csyn{\fy}&\leq\n{\fy+d(p_0\fy_*)}+\n{\fy+d(p_1\fy_*)}\\
&=\n{\fy_0+p_0(d\fy_*)}+\n{\fy_1}+\n{\fy_0}+\n{\fy_1+p_1(d\fy_*)}.
\end{split}
\end{equation}
Let $A$ denote the set of cells of $Q_{d-1}$ obtained from $\supp\fy_0$
by deleting the first coordinate (which is $0$), and let $B$ denote the
set of cells of $Q_{d-1}$ obtained from $\supp\fy_1$ by deleting the
first coordinate (which is $1$). Let $X$ denote the set of cells of
$Q_{d-1}$ obtained from $\supp(p_0(d\fy_*))$ by deleting the first
coordinate; note that it is the same set as the one obtained by
deleting the first coordinate in the elements of $\supp(p_1(d\fy_*))$.

The inequality~\eqref{eq:q1} says in this notation that
\[2\cdot\csyn{\fy}\leq\sn{A}+\sn{B}+\sn{A\xor X}+\sn{B\xor X}.\]
Applying the inequality~\eqref{eq:setin}, and cancelling out
the factor~$2$, we obtain
\begin{equation}
\label{eq:q2}
\csyn{\fy}\leq\sn{A}+\sn{B}+\sn{A\xor B\xor X}.
\end{equation}
On the other hand, we have
\begin{equation}\label{eq:q5}
\n{d\fy}=\n{d^0\fy_0}+\n{d^1\fy_1}+\n{\fy_0^*+\fy_1^*+d\fy_*},
\end{equation}
where $\fy_0^*$ is obtained from $\supp\fy_0$ by changing the first
coordinate to~$*$, and $\fy_1^*$ is obtained the same way. By
induction assumption for $(d-1,k)$, we have
$\sn{A}=\n{\fy_0}\leq\n{d^0\fy_0}$ and
$\sn{B}=\n{\fy_1}\leq\n{d^1\fy_1}$.  Furthermore, deleting the first
coordinate, we obtain
\[\n{\fy_0^*+\fy_1^*+d\fy_*}=\sn{A\xor B\xor X}.\]
Combining these with~\eqref{eq:q2} and~\eqref{eq:q5}, we obtain
\[\csyn{\fy}\leq \n{d\fy},\]
which finishes the proof.
\qed



\subsection{Expansion after taking product with a simplex}
\label{sub:cube}

Let $X$ be a cell complex and let $\smp{n-1}$ denote the
$(n-1)$-simplex on the vertex set $V=[n]=\{1,\dots,n\}$, $n\geq 2$. The
product complex $Y=X\times\smp{n-1}$ is again a~cell complex, whose
cells are products of the cells of $X$ with the simplices of
$\smp{n-1}$. Gromov proved that for all $k\geq 0$ we have
\[
h^k(X\times\smp{n-1}) \geq\min\bigg\{h^k(X),\frac{n-k-1}{k+2}\bigg\},
\]
see section~2.11 in \cite{Gr}.  Here we provide an
elementary proof of the following slightly stronger bound.
\begin{thm}\label{t:prod}
Let $X$ be a cell complex, and let $n\geq 2$. Then for all $k\geq 0$
we have the inequality
\begin{equation*}
h^k(X\times\smp{n-1})\geq\min\bigg\{h^k(X),\max\Big\{1,\frac{n}{k+2}\Big\}\bigg\}.
\end{equation*}
\end{thm}
\pr We start with some preliminary observations. Let $v$ be a
vertex in $V$.  For $1 \leq k \leq n-1$ let $T_v:C^k(Y) \rightarrow
C^{k-1}(Y)$ be the linear map defined as follows: For $\phi \in
C^k(Y)$ and a $(k-1)$-dimensional cell $\alpha \times \beta \in X(i)
\times \dn(j)$ where $i+j=k-1$ let
$$
T_v \phi(\alpha \times \beta)=\phi (\alpha \times [v, \beta]).
$$
\begin{claim}
\label{c:acom}
Let $\phi \in C^k(Y)$. Let $i+j=k$ and let
$\sigma =\alpha \times \beta \in X(i) \times \dn(j) \subset Y(k)$. Then:
\begin{equation}
\label{e:acom}
d_{k-1}T_v \phi (\sigma)+T_v d_k \phi (\sigma)=
\left\{
\begin{array}{ll}
\phi(\sigma) & (i,j) \neq (k,0), \\
\phi(\sigma)+\phi(\alpha \times v) & (i,j)=(k,0).
\end{array}
\right.~~
\end{equation}
\end{claim}
\pr
\begin{equation}
\label{e:left}
\begin{split}
d_{k-1}T_v \phi(\alpha \times \beta) &= T_v \phi\big( \partial_k(\alpha \times \beta)\big) \\
&= T_v \phi (\partial_i\alpha\times \beta+\alpha \times \partial_j \beta) \\
&=
\left\{
\begin{array}{ll}
\phi(\partial_i \alpha \times [v,\beta])+\phi(\alpha \times [v,\partial_j\beta])
& (i,j) \neq (k,0), \\
\phi(\partial_k \alpha \times [v,\beta]) & (i,j)=(k,0).
\end{array}
\right.~~
\end{split}
\end{equation}

\begin{equation}
\label{e:right}
\begin{split}
T_v d_k \phi(\alpha \times \beta) &=
d_k \phi(\alpha \times [v,\beta]) \\
&=\phi\big(\partial_{k+1}(\alpha \times [v,\beta])\big)
\\
&= \phi(\partial_i\alpha \times [v,\beta])+\phi(\alpha \times \partial_{j+1}[v,\beta]) \\
&=
\left\{
\begin{array}{ll}
\phi(\partial_i\alpha \times [v,\beta])+ \phi(\alpha\times \beta)+
\phi(\alpha\times [v,\partial_j\beta])  & (i,j) \neq (k,0), \\
\phi(\partial_k\alpha \times [v,\beta])+ \phi(\alpha\times \beta)+
\phi(\alpha \times v) & (i,j)=(k,0).
\end{array}
\right.~~
\end{split}
\end{equation}
Now (\ref{e:acom}) follows from (\ref{e:left}) and (\ref{e:right}).
{\enp}
\noindent
Claim \ref{c:acom} implies the following. If $\sigma =\alpha \times
\beta \in X(i) \times \dn(j)$ where $i+j=k$ and $0<j \leq \mkn$ then:
\begin{equation}
\label{e:jpos}
\begin{split}
(\phi+d_{k-1}T_v \phi)(\sigma) &= (\phi+d_{k-1}T_v \phi)(\alpha \times \beta) \\
&=T_vd_k\phi(\alpha \times \beta) \\
&=d_k \phi(\alpha\times [v,\beta]).
\end{split}
\end{equation}
\noindent
On the other hand, if $\sigma=\alpha \times u \in X(k) \times \dn(0)$ then:
\begin{equation}
\label{e:jzer}
\begin{split}
(\phi+d_{k-1}T_v \phi)(\sigma) &= (\phi+d_{k-1}T_v \phi)(\alpha \times u) \\
&=T_vd_k\phi(\alpha \times u)+\phi(\alpha \times v) \\
&=d_k \phi(\alpha\times [v,u])+ \phi(\alpha \times v).
\end{split}
\end{equation}
\noindent
For $0 \leq j \leq \mkn$ let
\[
g_j(\phi)=\sum_{v \in V}
\mybar{\supp(\phi+d_{k-1}T_v\phi)\cap\big(X(k-j)\times\dn(j)\big)}.
\]
By (\ref{e:jpos}) it follows that for every $1 \leq j \leq \mkn$
\begin{equation}
\label{e:sjpos}
g_j(\phi)=(j+2) \cdot\mybar{\supp(d_k\phi) \cap \big(X(k-j) \times \dn(j+1)\big)}.
\end{equation}
For $v \in V$ define the restriction map $R_v:C^k(Y) \rightarrow C^k(X)$ as follows.
For $\phi \in C^k(Y)$ and $\alpha \in X(k)$ let
$R_v\phi(\alpha)=\phi(\alpha \times v)$.
By (\ref{e:jzer}) it follows that
\begin{equation}
\label{e:sjzer}
\begin{split}
g_0&(\phi)=\sum_{v \in V}\mybar{\supp(\phi+d_{k-1}T_v \phi) \cap \big(X(k) \times \dn(0)\big)} \\
&=\mybar{\{(v,\alpha,u) \in \dn(0) \times X(k) \times \dn(0):d_k\phi(\alpha \times [v,u])
\neq \phi(\alpha \times v) \}} \\
&\leq\mybar{\mycb{(v,\alpha,u)\in\dn(0)\times X(k)\times\dn(0):d_k\phi(\alpha\times [v,u])\neq 0}} \\
&+ |\{(v,\alpha,u) \in \dn(0) \times X(k) \times \dn(0):R_v\phi(\alpha) \neq 0 \}| \\
&= 2\cdot\mybar{\supp(d_k \phi) \cap \big(X(k)\times \dn(1)\big)}+n\sum_{v \in V}\n{R_v\phi}.
\end{split}
\end{equation}
For $0 \leq \ell \leq m $ let
\[\cf_{m,\ell}=\bigcup_{j \geq \ell} X(m-j) \times \dn(j).\]
Combining (\ref{e:sjpos}) and (\ref{e:sjzer}) we obtain
\begin{equation}
\label{e:ubphi}
\begin{split}
n&\csyn{\phi} \leq \sum_{v \in V} \n{\phi+d_{k-1}T_v \phi} \\
&=\sum_{v \in V} \sum_{j=0}^{\mkn}\mybar{\supp(\phi+d_{k-1}T_v \phi) \cap \big(X(k-j)
\times \dn(j)\big)} \\
&= \sum_{j=0}^{\mkn} g_j(\phi) \\
&\leq \sum_{j=0}^{\mkn} (j+2) \cdot\mybar{\supp(d_k\phi) \cap \big(X(k-j)
\times \dn(j+1)\big)}+n\sum_{v \in V}\|R_v\phi\| \\
&\leq \min\{k+2,n\}\cdot\mybar{\supp(d_k\phi) \bigcap \cf_{k+1,1}}+
n \sum_{v \in V} \|R_v\phi\|.
\end{split}
\end{equation}
\begin{claim}
\label{c:local}
For any $\phi \in C^k(Y)$, there exists a $\psi \in C^{k-1}(Y)$, such
that $\tphi=\phi+d_{k-1} \psi$ satisfies $\n{R_v
  \tphi}=\csyn{R_v\phi}$ for all $v \in V$.
\end{claim}
\noindent
{\bf Proof:} For $v \in V$ choose $\psi_v \in C^{k-1}(X)$ such that
$\|R_v \phi\|_{csy}=\|R_v\phi+d_{k-1} \psi_v\|$.  Define $\psi \in
C^{k-1}(Y)$ by
$$
\psi(\sigma)=
\left\{
\begin{array}{ll}
\psi_v(\alpha) & \sigma=\alpha \times v \in X(k-1) \times \dn(0), \\
0 & otherwise.
\end{array}
\right.~~
$$
Noting that $R_v \psi=\psi_v$ and $R_v d_{k-1}=d_{k-1}R_v$, it follows that
\begin{equation*}
\label{e:tilde}
\begin{split}
\n{R_v \tphi} &=\|R_v\phi+R_v d_{k-1}\psi\| \\
&=\|R_v \phi+d_{k-1}R_v \psi\| \\
&=\|R_v\phi+d_{k-1}\psi_v\| \\
&= \|R_v \phi\|_{csy}.
\end{split}
\end{equation*}
{\enp}
Claim \ref{c:local} implies that
\begin{equation}
\label{e:locc}
\begin{split}
\sum_{v \in V}\n{R_v\tphi} &=\sum_{v \in V} \|R_v \phi\|_{csy} \\
&\leq \frac{1}{h^k(X)}\sum_{v \in V} \n{d_k R_v \phi} \\
&= \frac{1}{h^k(X)}\mybar{\supp(d_k \phi) \cap \big(X(k+1)\times \dn(0)\big)}.
\end{split}
\end{equation}
Applying (\ref{e:ubphi}) for $\tphi$ and using (\ref{e:locc}) we obtain:
\begin{equation}
\label{e:finbd}
\begin{split}
n\|\phi\|_{csy}&=n\csyn{\tphi} \\
&\leq \min\{k+2,n\}\cdot\mybar{\supp(d_k\tphi)\bigcap\cf_{k+1,1}}+
n \sum_{v \in V}\n{R_v\tphi} \\
&\leq \min\{k+2,n\}\cdot\mybar{\supp(d_k\phi) \bigcap \cf_{k+1,1}} \\
&~~~~~~~+\frac{n}{h^k(X)}\cdot\mybar{\supp(d_k \phi) \cap \big(X(k+1)\times \dn(0)\big)} \\
&\leq \max\mycb{\min\{k+2,n\},\frac{n}{h^k(X)}}\cdot \|d_k \phi\|.
\end{split}
\end{equation}
\noindent
Rearranging (\ref{e:finbd}) it follows that
\[
\frac{\|d_k \phi\|}{\|\phi\|_{csy}} \geq \min
\mycb{h^k(X),\max\mycb{1,\frac{n}{k+2}}}.
\]
{\enp}

Note that Theorem~\eqref{thm:cube} is a direct corollary
of Theorem~\ref{t:prod}.
\editcut

\section{Maximal $k$-Cosystoles and Maximal Cheeger Constants}
\label{s:maxs}
\mysubs{Systolic norm and expansion of random cochains}
\noindent
Let $X$ be a simplicial complex and let $0 \leq k \leq \dim X$.
Let $\lambda_k(X)$ denote the maximal norm of a $k$-cosystole in $X$:
\[\lambda_k(X)=\max\mycb{\|\phi\|_{csy}: \phi \in C^k(X)}.\]
\begin{theorem}
\label{p:lkn}
\begin{equation}
\label{e:blam}
\Bigg(1-20\sqrt{\frac{f_{k-1}(X)}{f_k(X)}}\Bigg)\cdot \frac{f_k(X)}{2}
\leq \lambda_k(n) \leq \frac{f_k(X)}{2}.
\end{equation}
\end{theorem}
\noindent
{\bf Proof:} The upper bound is straightforward.
Let $\phi \in C^k(X)$ be a $k$-cosystole and let $\tau \in X(k-1)$.
Write $\Gamma_X(\tau)=\{\sigma \in X(k): \sigma \supset \tau\}$ and recall that
$\deg_X(\tau)=|\Gamma_X(\tau)|$. It follows that
\begin{equation}
\label{e:lthf}
0 \leq\n{\phi+d_{k-1}\tau^*}-\n{\phi}= \deg_X(\tau)- 2|\Gamma_X(\tau) \cap \supp(\phi)|.
\end{equation}
Rearranging and summing (\ref{e:lthf}) over all $\tau \in X(k-1)$ we obtain
\begin{equation*}
\label{e:lbnd}
\begin{split}
(k+1)\|\phi\|&= \sum_{\tau \in X(k-1)}|\Gamma_X(\tau) \cap \supp(\phi)| \\
&\leq \frac{1}{2}\sum_{\tau \in X(k-1)} \deg_X(\tau)=\frac{1}{2}(k+1)f_k(X).
\end{split}
\end{equation*}

We next prove the lower bound. Let $N=f_k(X), M=f_{k-1}(X)$ and let
$\delta=10\sqrt{\frac{M}{N}}$.
Consider the probability space of all $k$-cochains
\[\phi=\sum_{\sigma \in X(k)} x_{\sigma} \sigma^* \in C^k(X),\]
where $\left\{x_{\sigma}:\sigma \in X(k)\right\}$ is a family of independent $0,1$ variables with
$\prob[x_{\sigma}=0]=\prob[x_{\sigma}=1]=\frac{1}{2}$.
\begin{claim}
\label{c:nearcsy}
\begin{equation}
\label{e:nearcsy}
\prob\left[\|\phi\|_{csy} < \left(1-20\sqrt{\frac{M}{N}}\right)\frac{N}{2} \right] \leq 0.8^M.
\end{equation}
\end{claim}
\noindent
{\bf Proof:} If $\delta=10\sqrt{\frac{M}{N}} \geq \frac{1}{2}$ then the claim is vacuous, so we shall henceforth assume that $\delta<\frac{1}{2}$. Let $Y$ be the random variable given
by $Y(\phi)=\|\phi\|$. Then $E[Y]=\frac{N}{2}$ and by Chernoff's bound (see, e.g., \cite{AS})
\begin{equation}
\label{e:cher1}
\prob\left[Y<\frac{N}{4}\right]< e^{-\frac{N}{8}}.
\end{equation}
Fix a $(k-1)$-cochain $\psi \in C^{k-1}(X)$ and let
$S(\psi)=\supp(d_{k-1} \psi) \subset X(k)$.  Let $Z_{\psi}$ be the
random variable given by
\[
Z_\psi(\phi)=\sn{S(\psi) \cap \supp(\phi)}=\sn{\{\sigma \in S(\psi): x_{\sigma}=1\}}.
\]
Then $E\left[Z_{\psi}\right]=\frac{|S(\psi)|}{2}$. By Chernoff's bound
\begin{equation}
\label{e:cher2}
\prob\left[Z_{\psi} \geq \frac{|S(\psi)|}{2}+\frac{\delta N}{8}\right] \leq
e^{-\frac{\delta^2 N^2}{32 |S(\psi)|}} \leq e^{-\frac{\delta^2 N}{32}}.
\end{equation}
Next note that
\begin{equation}
\label{e:mcc}
\begin{split}
\|\phi+d_{k-1}\psi\|&=\|\phi\|+|S(\psi)|-2|S_{\psi} \cap \supp(\phi)| \\
&=Y(\phi)+|S(\psi)|-2 Z_{\psi}(\phi).
\end{split}
\end{equation}
Combining (\ref{e:mcc}), (\ref{e:cher1}) and (\ref{e:cher2}) we obtain
\begin{equation}
\label{e:prsm}
\begin{split}
\prob\left[\n{\phi+d_{k-1}\psi} \leq (1-\delta)\n{\phi}\right]
&=\prob\left[Y+|S(\psi)|-2 Z_{\psi} \leq (1-\delta)Y\right] \\
&=\prob\left[Z_{\psi} \geq \frac{|S(\psi)|}{2} +\frac{\delta Y}{2}\right] \\
&\leq\prob\left[Y \leq \frac{N}{4}\right]+\prob\left[Z_{\psi} \geq
\frac{|S(\psi)|}{2}+\frac{\delta N}{8}\right] \\
&\leq e^{-\frac{N}{8}}+e^{-\frac{\delta^2 N}{32}} \leq 2e^{-\frac{\delta^2 N}{32}}.
\end{split}
\end{equation}
Therefore
\begin{equation}
\label{e:psmall}
\begin{split}
&\prob\left[\|\phi\|_{csy} \leq (1-2\delta)\frac{N}{2}\right] \\
&\leq \prob\left[\|\phi\|\leq (1-\delta)\frac{N}{2}\right]+
\sum_{\psi \in C^{k-1}(X)}\prob\left[\|\phi+d_{k-1}\psi\| \leq (1-\delta)\|\phi\|\right] \\
&\leq e^{-\frac{\delta^2 N}{2}}+2^M \cdot 2e^{-\frac{\delta^2 N}{32}} \\
&=e^{-50M}+2^{M+1}e^{-\frac{100M}{32}} <0.8^M.
\end{split}
\end{equation}
{\enp}
Claim \ref{c:nearcsy} implies that there exists a $\phi \in C^k(X)$ such that
\[
\|\phi\|_{csy} \geq (1-2\delta)\frac{N}{2}=\Bigg(1-20\sqrt{\frac{f_{k-1}(X)}{f_k(X)}}\Bigg)\cdot \frac{f_k(X)}{2}.
\]
{\enp}
\noindent
 Claim \ref{c:nearcsy} can also be used to provide an upper bound on the $k$-th Cheeger constant of certain sparse complexes.
\begin{theorem}
\label{t:uphk}
Let $X$ be a pure $(k+1)$-dimensional complex such that $\deg_X(\sigma)=D$ for every $\sigma \in X(k)$.
If $D\geq 40^2 (k+1)$ then
\begin{equation}
\label{e:uphk}
h^k(X) \leq \left(1+\frac{50 \sqrt{k+1}}{\sqrt{D}}\right) \cdot \frac{D}{k+2}.
\end{equation}
\end{theorem}
For the proof we will need the following consequence of Azuma's inequality due to McDiarmid \cite{McD}.
\begin{theorem}
\label{t:azuma}
Suppose $g:\{0,1\}^N \rightarrow \Rea$ satisfies
$|g(\epsilon)-g(\epsilon')| \leq D$ if $\epsilon$ and $\epsilon'$ differ in at most one coordinate.
Let $x_1,\ldots,x_N$ be independent $0,1$ valued random variables
and let $G=g(x_1,\ldots,x_N)$. Then for all $\lambda>0$
\begin{equation}
\label{azumain}
\prob[G \leq E[G] - \lambda] \leq \exp\left(-\frac{2\lambda^2}{D^2 N}\right)~.
\end{equation}
\end{theorem}
{\bf Proof of Theorem \ref{t:uphk}:}
Let $L=f_{k+1}(X), N=f_k(X), M=f_{k-1}(X)$.
Let $X(k)=\{\sigma_1,\ldots,\sigma_N\}$ and let $x_1,\ldots,x_N$ be
independent $0,1$ variables with
$\prob[x_i=0]=\prob[x_i=1]=\frac{1}{2}$.
Define $g:\{0,1\}^N \rightarrow \Rea$ by
$$g(\epsilon_1,\ldots,\epsilon_N)=L-\|d_k\left(\sum_{i=1}^N \epsilon_i \sigma_i^*\right)\|.$$
The assumption $\deg_X(\sigma)=D$ for all $\sigma \in X(k)$ implies that $|g(\epsilon)-g(\epsilon')| \leq D$ if $\epsilon$ and $\epsilon'$ differ in at most one coordinate. Let $\phi$ be the random $k$-cochain $\phi=\sum_{i=1}^N x_i \sigma_i^*$.
Then $G=g(x_1,\ldots,x_N)$ satisfies $E[G]=L-E[\|d_k \phi\|]=\frac{L}{2}$.
Thus, by Theorem \ref{t:azuma}
\begin{equation}
\label{e:ubdkg}
\begin{split}
\prob\left[\|d_k\phi\| \geq \frac{L}{2}+D\sqrt{N}\right]
&=\prob\left[G \leq \frac{L}{2}-D\sqrt{N}\right] \\
&< \exp\left(-\frac{2(D\sqrt{N})^2}{D^2N}\right)=\exp(-2)<0.2.
\end{split}
\end{equation}
Combining (\ref{e:ubdkg}) and (\ref{e:nearcsy}) it follows that there exists a $\phi \in C^k(X)$ such that
\begin{equation}
\label{e:nomi}
\|d_k\phi\| \leq \frac{L}{2}+D\sqrt{N}
\end{equation}
and
\begin{equation}
\label{e:denomi}
\|\phi\|_{csy} \geq \left(1-20\sqrt{\frac{M}{N}}\right)\cdot \frac{N}{2}.
\end{equation}
Next note that
\begin{equation}
\label{e:eas1}
(k+2)L=DN
\end{equation}
and
\begin{equation}
\label{e:eas2}
N \geq 1+D(k+1).
\end{equation}
Furthermore, $\deg_X(\tau) \geq D+1$ for all $\tau \in X(k-1)$
and therefore
\begin{equation}
\label{e:eas3}
(k+1)N \geq (D+1)M.
\end{equation}
Combining (\ref{e:nomi}) and (\ref{e:denomi}), and using
(\ref{e:eas1}), (\ref{e:eas2}), (\ref{e:eas3}) and the assumption $\sqrt{\frac{k+1}{D}} \leq \frac{1}{40}$, we obtain
\begin{equation}
\label{e:ubhkk}
\begin{split}
h^k(X) &\leq \frac{\|d_k\phi\|}{\|\phi\|_{csy}}
\leq \frac{\frac{L}{2}+D\sqrt{N}}{\left(1-20\sqrt{\frac{M}{N}}\right) \cdot \frac{N}{2}} \\
&= \frac{\left(1+\frac{2(k+2)}{\sqrt{N}}\right)\cdot\frac{L}{2}}{\left(1-20\sqrt{\frac{M}{N}}\right) \cdot \frac{N}{2}} \\
&\leq \frac{1+\frac{2(k+2)}{\sqrt{N}}}{1-20\sqrt{\frac{k+1}{D+1}}}\cdot \frac{D}{k+2} \\
&\leq \frac{1+\frac{4(k+1)}{\sqrt{D(k+1)}}}{1-20\sqrt{\frac{k+1}{D+1}}}\cdot \frac{D}{k+2} \\
&\leq \left(1+4\sqrt{\frac{k+1}{D}}\right)\left(1+40\sqrt{\frac{k+1}{D+1}}\right)\cdot \frac{D}{k+2} \\
&\leq \left(1+50\sqrt{\frac{k+1}{D}}\right)\cdot \frac{D}{k+2}
\end{split}
\end{equation}
{\enp}

\mysubs{The cosystolic norm of the Paley cochain}

Let $p>2$ be a~prime and let $\chi$ be the quadratic character of $\FF_p$,
i.e., $\chi(x)=\myrb{\frac{x}{p}}$, the Legendre symbol of $x$ modulo $p$.
The \emph{Paley graph} $G_p$ is the graph on the vertex set $\FF_p$
whose edges are pairs $\{x,y\}$ such that $\chi(x-y)=1$. The Paley graph is an important example of an explicitly given graph that exhibits strong pseudorandom properties (see, e.g., \cite{AS}).
Motivated by the above, we now define a high dimensional version of the Paley graph.
Let $1 \leq k < p$ be fixed and let $\smp{p-1}$ be the $(p-1)$-simplex on the
vertex set $\FF_p$.
The {\bf Paley $k$-Cochain} $\phi_k \in C^k(\smp{p-1})$ is defined as follows. For a
$k$-simplex $\sigma=\{x_0,\ldots,x_k\}$ let $\phi_k(\sigma)=1$ if
$\chi\myrb{x_0+\dots+x_k}=1$, and $\phi_k(\sigma)=0$ otherwise.
Here we prove that the Paley $k$-cochain $\phi_k$ is close to being a $k$-cosystole in $C^k(\smp{p-1})$.
\begin{theorem}
\label{p:paley}
For a fixed $k \geq 1$
$$\|\phi_k\|_{csy} \geq \frac{1}{2} \binom{p}{k+1} \myrb{1-O\myrb{p^{-2^{-k}}}}.$$
\end{theorem}
\noindent
The proof of Theorem \ref{p:paley} depends on the following result
of Chung \cite{Ch}.  For $0 \leq i \leq k$ define the projection
$\pi_i:\FF_p^{k+1} \rightarrow \FF_p^k$ by
$\pi_i(x_0,\ldots,x_k)=(x_0,\ldots,x_{i-1},x_{i+1},\ldots,x_k)$. For
subsets $R_0,\ldots,R_k \subset \FF_p^k$ let
\[W(R_0,\ldots,R_k)=\mycb{x \in \FF_p^{k+1}: \pi_i(x) \in R_i \text{~for~all~} 0 \leq i \leq k}.\]
\begin{thm}[Chung \cite{Ch}]
\label{t:chung}
\begin{equation}
\label{e:chung}
\begin{split}
&\Big|\sum_{\sigma=(x_0,\ldots,x_k) \in W(R_0,\ldots,R_k)}
\chi(x_0+\dots+x_k)\Big| \\
&\leq 2^{(k-1)2^{-(k-1)}} p^{1-2^{-k}}\myrb{\prod_{i=0}^k |R_i|}^{\frac{1}{k+1}}
\leq 2 p^{k+1-2^{-k}}.
\end{split}
\end{equation}
\end{thm}
\prn{Theorem \ref{p:paley}}  First note that sums of $k+1$ distinct
elements of $\FF_p$ are equidistributed in $\FF_p$, hence
\begin{equation}
\label{e:phinorm}
\begin{split}
\|\phi_k\|&=\big|\big\{\{x_0,\ldots,x_k\} \in \ddp(k): \chi\myrb{x_0+\dots+x_k}=1\big\}\big| \\
&=\frac{\big|\big\{y \in \FF_p: \chi(y)=1\big\}\big|}{p}\binom{p}{k+1}=\frac{p-1}{2p} \binom{p}{k+1}.
\end{split}
\end{equation}
Let $\psi\in C^{k-1}(\smp{p-1})$ such that
$\|\phi\|_{csy}=\|\phi+d_{k-1}\psi\|$.
Let $S(\psi)=\supp(d_{k-1} \psi) \subset \dn(k)$.
Then
\begin{equation}
\label{e:diffr}
\begin{split}
&\sn{S(\psi)}-2\sn{\supp(\phi)\cap S(\psi)} \\
&=\sn{S(\psi)\setminus\supp(\phi)}-\sn{S(\psi)\cap\supp(\phi)} \\
&=\sn{\big\{\{x_0,\ldots,x_k\}\in S(\psi):\chi\myrb{x_0+\dots+x_k}\neq 1\big\}}\\
&\,\,\,\,\,-\sn{\big\{\{x_0,\ldots,x_k\}\in S(\psi):\chi\myrb{x_0+\dots+x_k}=1\big\}} \\
&=\sn{\big\{\{x_0,\ldots,x_k\} \in S(\psi): x_0+\dots+x_k=0\big\}} \\
&\,\,\,\,-\sum_{\{x_0,\ldots,x_k\}\in S(\psi)}\!\!\!\!\! \chi\myrb{x_0+\dots+x_k}\\
&\geq -\sum_{\{x_0,\ldots,x_k\} \in S(\psi)}\!\!\!\!\! \chi\myrb{x_0+\dots+x_k}.
\end{split}
\end{equation}
We proceed to bound the sum
\[\big|\sum_{\{x_0,\ldots,x_k\} \in S(\psi)} \chi\myrb{x_0+\dots+x_k}\big|.\]
Let
\[A=\mycb{(y_1,\ldots,y_k) \in \FF_p^k:\{y_1,\ldots,y_k\} \in \supp(\psi)}\]
and let
\[\tilde{S}(\psi)=\mycb{(x_0,\ldots,x_k) \in \FF_p^{k+1}:
\{x_0,\ldots,x_k\} \in S(\psi)}.\]
Denote
\[D=\mycb{(y_1,\ldots,y_k)\in\FF_p^k: y_i \neq y_j \text{~for~all~} i \neq j}\]
and let $$A_0=D\setminus A\,\,\,\, ,\,\,\,\, A_1=A.$$  For $\ueps=(\eps_0,\ldots,\eps_k) \in
\{0,1\}^{k+1}$ let $$W(\ueps)=W(A_{\eps_0},\ldots,A_{\eps_k}).$$
Write
\[E=\Big\{\ueps=(\eps_0,\ldots,\eps_k) \in \{0,1\}^{k+1}:
\eps_0+\dots+\eps_k \equiv 1(\text{mod} ~2)\Big\}.\]
Then
\begin{equation}
\label{e:base}
\tilde{S}(\psi)=\bigcup_{\ueps \in E} W(\ueps).
\end{equation}
By \eqref{e:base} and \eqref{e:chung} we have
\begin{equation}
\label{e:esti}
\begin{split}
&\Big|\sum_{\{x_0,\ldots,x_k\} \in S(\psi)} \chi\myrb{x_0+\dots+x_k}\Big| \\
&=\frac{1}{(k+1)!}\Big|\sum_{(x_0,\ldots,x_k) \in \tilde{S}(\psi)} \chi\myrb{x_0+\dots+x_k}\Big|\\
&\leq \frac{1}{(k+1)!} \sum_{\ueps \in E}\Big|\sum_{(x_0,\ldots,x_k) \in W(\ueps)} \chi\myrb{x_0+\dots+x_k}\Big| \\
&\leq \frac{2^{k+1}}{(k+1)!} p^{k+1-2^{-k}}.
\end{split}
\end{equation}
Combining (\ref{e:diffr}), (\ref{e:phinorm}) and (\ref{e:esti}), we obtain
\begin{equation*}
\label{e:estf}
\begin{split}
\|\phi_k\|_{csy}&=\|\phi_k+d_{k-1}\psi\| \\
&= \|\phi_k\|+|S(\psi)|-2|S(\psi) \cap \supp(\phi_k)| \\
&\geq\|\phi_k\|-\bigg|\sum_{\{x_0,\ldots,x_k\}\in S(\psi)}\chi\myrb{x_0+\dots+x_k}\bigg|\\
&\geq\frac{p-1}{2p}\binom{p}{k+1}-\frac{2^{k+1}}{(k+1)!} p^{k+1-2^{-k}} \\
&\geq \frac{1}{2} \binom{p}{k+1} \myrb{1-O\myrb{p^{-2^{-k}}}}. \qed
\end{split}
\end{equation*}

\noindent
{\bf Remark:} In the graphical case $k=1$, the Paley $1$-cochain
$\phi_1$ satisfies
$\|\phi_1\|_{csy}=\frac{1}{2}\binom{p}{2}\myrb{1-O\myrb{p^{-\frac{1}{2}}}}$.
It would be interesting  to decide whether $\|\phi_k\|_{csy} = \frac{1}{2}
\binom{p}{k+1} \myrb{1-O\myrb{p^{-\frac{1}{2}}}}$ remains true for
$k\geq 2$ as well.

\section{Bounded Quotients of the Fundamental Group of a Random $2$-Complex}
\label{s:nonab}

\mysubs{Probability space of simplicial complexes}

Let $Y(n,p)$ denote the probability space of random $2$-dimensional subcomplexes of $\dn$ obtained by starting
with the full $1$-skeleton of $\dn$ and then adding
each $2$-simplex independently with probability $p$.
Formally, $Y(n,p)$ consists of all
complexes $(\dn)^{(1)} \subset Y \subset (\dn)^{(2)}$ with
probability measure
$$\Pr(Y)=p^{f_2(Y)}(1-p)^{\binom{n}{3}-f_2(Y)}~.$$

Note that $p$ is a function of $n$, which is typically not
a~constant. Still, to simplify notations we omit the argument, and
just write $p$ instead of $p(n)$.

The threshold probability
for the vanishing of the first homology with fixed finite abelian
coefficient group $R$ was determined in \cite{LiM,MW}.
\begin{theorem}[\cite{LiM,MW}]
\label{gen}  Let $R$ be an~arbitrary finite abelian group, and let
$\omega(n):{\mathbb N}\ra {\mathbb R}$ be an arbitrary function that
goes to infinity when $n$ goes to infinity. Then the following
asymptotic result holds
\[
\lim_{n \rightarrow \infty} \prob ~[~X \in Y(n,p):
H_1(X;R)=0~]=
\begin{cases}
        0, & \textrm{ if } p=\frac{2\log n -\omega(n)}{n}; \\
        1, & \textrm{ if } p=\frac{2\log n+ \omega(n)}{n}.
\end{cases}
\]
\end{theorem}
The case of integral homology was addressed by Hoffman, Kahle and Paquette \cite{HKP} who proved that there exists a constant $c$ such that if
$p>\frac{c \log n}{n}$ then $X \in Y(n,p)$ satisfies $H_1(X;\Int)=0$ asymptotically almost surely.
Recently, {\L}uczak and Peled \cite{LP} proved that $p=\frac{2 \log n}{n}$ is a sharp threshold for the vanishing of $H_1(X;\Int)$.

Similarly, one can ask what is the threshold probability for the
vanishing of the fundamental group in the probability space $Y(n,p)$.
This is a quite difficult question which was answered by Babson,
Hoffman and Kahle, see~\cite{BHK}.
\begin{theorem}[\cite{BHK}]
\label{bhkt}
Let $\varepsilon>0$ be fixed, then
\[
\lim_{n \rightarrow \infty} \prob ~[~X \in Y(n,p):
\pi_1(X)=0 ~]=
\begin{cases}
        0, & \textrm{ if } p=\left(\frac{n^{-\varepsilon}}{n}\right)^{1/2}; \\
        1, & \textrm{ if } p=\left(\frac{3\log n +\omega(n)}{n}\right)^{1/2}.
\end{cases}
\]
\end{theorem}
In view of the gap between the thresholds for the vanishing of $H_1(X;\Int)$ and 
for the triviality of $\pi_1(Y)$, Eric Babson (see problem (8) on page~58
in \cite{FGS}) asked what is the threshold probability such that a.a.s.\ $\pi_1(X)$
does not have a quotient equal to some non-trivial finite group. This a property
between the vanishing of the fundamental group and the vanishing of the first
homology group. If the fundamental group is trivial, then certainly it cannot
contain a~quotient equal to a~non-trivial finite group. On the other hand, if
the first homology group is non-trivial, then it must be finite, and it is
the quotient of the fundamental group by its commutator, so this property is satisfied.

Addressing Babson's question we prove the following theorem.
\begin{theorem}
\label{nonab}
Assume $c>0$ is a~constant, and set $p:=\frac{(6+7c)\log n}{n}$.
Then, when $X$ is sampled from the probability space $Y(p,n)$,
a.a.s.\ the fundamental group $\pi_1(X)$ does not contain a~proper
normal subgroup of index at most~$n^c$.
\end{theorem}

{\bf Remark.} The constant $6+7c$ may be improved using a more careful
analysis as in \cite{LiM,MW}. For example, for any {\it fixed}
non-trivial finite group $G$, if $p=\frac{2\log n+\omega(n)}{n}$ then
a.a.s.\ $G$ is not a homomorphic image of $\pi_1(X)$.

The proof of Theorem \ref{nonab} is an adaptation of the argument in
\cite{LiM,MW} to the non-abelian setting. In Section \ref{naco} we
recall the notion of non-abelian first cohomology and its relation
with the fundamental group.  In Section \ref{expsim} we compute the
expansion of the $(n-1)$-simplex.  The results of Sections \ref{naco}
and \ref{expsim} are used in section \ref{fmw} to prove Theorem
\ref{nonab}.

\mysubs{Non-abelian first cohomology}
\label{naco}
Let $X$ be a simplicial complex and let $G$ be a multiplicative group.
We do not assume that $G$ is abelian. The definition of the first
cohomology $H^1(X;G)$ of $X$ with coefficients in $G$ was given
in~\cite{Olum58}. Since the setting in \cite{Olum58} was that of
singular cohomology, whereas we would like to work simplicially, we
choose to include a~certain amount of details in our recollection
below.

For $0 \leq k \leq 2$, let $\wti X(k)$ denote the set of all ordered
$k$-simplices of $X$.  Let $C^0(X;G)$ denote the group of $G$-valued
functions on $\wti X(0)=X(0)$ with pointwise
multiplication. Furthermore, set
\[
C^1(X;G):=\{\phi:\wti X(1) \rightarrow G: \phi(u,v)=\phi(v,u)^{-1}\}.
\]
We define the $0$-th coboundary operator $d_0:C^0(X;G)\rightarrow
C^1(X;G)$ by setting
\[
(d_0\psi)(u,v):=\psi(u)\,\psi(v)^{-1},
\]
for all $\psi\in C^0(X,G)$, and $(u,v)\in\wti X(1)$.

Proceeding to dimension $2$, let
$C^2(X;G)$ denote the set of all functions $\{\phi:\wti X(2) \rightarrow G$.
Define the first
coboundary operator $d_1:C^1(X;G)\rightarrow C^2(X;G)$ by setting
\[
(d_1\phi)(u,v,w):=\phi(u,v)\,\phi(v,w)\,\phi(w,u),
\]
for all $\phi\in C^1(X;G)$ and $(u,v,w)\in\wti X(2)$.

Define the set of $G$-valued $1$-cocycles of $X$ by
\[
Z^1(X;G):=\{\phi \in C^1(X;G):(d_1\phi)(u,v,w)=1,{\rm ~for~all~}(u,v,w) \in\wti X(2)\}.
\]
Furthermore, define an action of $C^0(X;G)$ on $C^1(X;G)$ by
setting
\begin{equation}
\label{eq:dotact}
(\psi\, .\,\phi)(u,v)=\psi(u)\,\phi(u,v)\,\psi(v)^{-1},
\end{equation}
for all $\psi \in C^0(X;G)$ and all $\phi \in C^1(X;G)$. In
particular, we recover the $0$-th coboundary operator as $d_0 \psi
=\psi\, .\, 1$. For $\phi \in C^1(X;G)$ let $[\phi]$ denote the orbit
of $\phi$ under that action.

We claim that $Z^1(X;G)$ is invariant under the action of $C^0(X;G)$. Indeed,
let $\phi\in Z^1(X;G)$, then for all $\psi\in C^0(X;G)$
\begin{equation*}
\begin{split}
d_1(\psi\,.\,\phi)(u,v,w)&=(\psi\,.\,\phi)(u,v)\ (\psi\,.\,\phi)(v,w)\
(\psi\,.\,\phi)(w,u)\\
&=\psi(u)\phi(u,v)\psi(v)^{-1}\,\psi(v)\phi(v,w)\psi(w)^{-1}\,
\psi(w)\phi(w,u)\psi(u)^{-1}\\
&=\psi(u)\,\phi(u,v)\,\phi(v,w)\,\phi(w,u)\,\psi(u)^{-1}\\
&=\psi(u)\, 1\,\psi(u)^{-1}=1.
\end{split}
\end{equation*}

We can now define non-abelian cohomology in dimension~$1$.

\begin{df} \label{df:nonabc}
The {\bf first non-abelian cohomology} of $X$ with coefficients in $G$
is the set of orbits of $Z^1(X;G)$ under the action of $C^0(X;G)$:
\[
H^1(X;G):=\{[\phi]: \phi \in Z^1(X;G) \}.
\]
\end{df}

Note that in general $H^1(X;G)$ is just a~set.
Furthermore, when $G$ is an abelian group, Definition~\ref{df:nonabc} yields
the usual first cohomology group of $X$ with coefficients in~$G$.

Assume now that the simplicial complex $X$ is connected, and let
$\Hom(\pi_1(X),G)$ denote the set of group homomorphisms from
$\pi_1(X)$ to $G$.  The group $G$ acts by conjugation
on $\Hom(\pi_1(X),G)$: for $\varphi \in \Hom(\pi_1(X),G)$ and
$g\in G$, let $g(\varphi) \in \Hom(\pi_1(X),G)$ be given
\[g(\varphi)(\gamma):=g\cdot\varphi(\gamma)\cdot g^{-1},\] for all $\gamma\in\pi_1(X)$.
For $\varphi \in \Hom(\pi_1(X),G)$ let $[\varphi]$ denote the orbit of
$\varphi$ under this action, and let
\[
\Hom(\pi_1(X),G)/G:=\{[\varphi]: \varphi \in \Hom(\pi_1(X),G)\}.
\]
The following observation is well known (see (1.3) in \cite{Olum58}).
For completeness we outline a proof.

\begin{prop}
\label{prop:fhpi}
For any  $(\dn)^{(1)} \subset X \subset (\dn)^{(2)}$
there exists a bijection
\[
\mu:\Hom(\pi_1(X),G)/G \rightarrow H^1(X;G),
\]
that maps $[1] \in \Hom(\pi_1(X),G)/G$ to $[1] \in H^1(X;G)$.
\end{prop}

\pr
We identify $\pi_1(X)$ with the group $\langle{\mathcal E}\mid{\mathcal R}\rangle$,
where the generating set is
\[{\mathcal E}:=\{e_{ij}: 2 \leq i, j \leq n,\, i\neq j\}\]
and set of relations ${\mathcal R}$ is given by
\begin{enumerate}
\item[(R1)]
$e_{ij}e_{ji}=1$, for all $i,j$,
\item[(R2)]
$e_{ij}=1$, if $(1,i,j) \in\wti X(2)$,
\item[(R3)]
$e_{ij}e_{jk}e_{ki}=1$, if $(i,j,k) \in\wti X(2)$.
\end{enumerate}
Each generator $e_{ij}$ corresponds to the loop consisting of $3$ edges:
$(1,i)$, $(i,j)$, and $(j,1)$, making the relations (R1)-(R3) obvious.

In these notations, each group homomorphism $\varphi:\pi_1(X)\ra G$ is induced by a set map
$\varphi:{\mathcal E}\ra G$ that maps all the relations (R1)-(R3) to the unit.
Furthermore, the conjugation action of $G$ on $\Hom(\pi_1(X),G)$ is induced by
$g(\varphi)(e_{ij})=g\,\varphi(e_{ij})\,g^{-1}$, for all $g\in G$, $\varphi:\pi_1(X)\ra G$.

For an~arbitrary group homomorphism $\varphi\in \Hom(\pi_1(X),G)$, define the cochain
$F(\varphi) \in C^1(X;G)$ by setting
\[
F(\varphi)(i,j):=
\begin{cases}
        \varphi(e_{ij}), &\textrm{ if } 2 \leq i, j \leq n,\, i\neq j, \\
        1, &\textrm{ otherwise, }
\end{cases}
\]
for all $1\leq i,j\leq n$, $i\neq j$. This is well-defined because
$F(\varphi)(i,j)=F(\varphi)(j,i)^{-1}$: if $i=1$ or $j=1$ this is trivial as both
sides are equal to~$1$,
and if $i,j\neq 1$ we get
\[F(\varphi)(j,i)=\varphi(e_{ji})=\varphi(e_{ij}^{-1})=\varphi(e_{ij})^{-1},\]
where the first equality is the definition of $F$, the second equality
follows from (R1), and the last equality follows from the fact that
$\varphi$ is a group homomorphism.

Let us see that $F(\varphi) \in Z^1(X;G)$, for all
$\varphi\in\Hom(\pi_1(X),G)$. Take $(u,v,w)\in\wti X(2)$. If $u=1$, then
\begin{equation*}
\begin{split}
d_1(F(\varphi))(1,v,w)&=F(\varphi)(1,v)\,F(\varphi)(v,w)\,F(\varphi)(w,1)\\
&=F(\varphi)(v,w)=\varphi(e_{v,w})=1,
\end{split}
\end{equation*}
where the last equality
follows from (R2). If $u\neq 1$, we can assume without loss of generality that also
$v,w\neq 1$. In that case we have
\begin{equation*}
\begin{split}
d_1(F(\varphi))(u,v,w)&=F(\varphi)(u,v)\,F(\varphi)(v,w)\,F(\varphi)(w,u)\\
&=\varphi(e_{u,v})\,\varphi(e_{v,w})\,\varphi(e_{w,u})=1,
\end{split}
\end{equation*}
where the last equality follows from (R3).

Let us see that the mapping
\[\wti{F}:\Hom(\pi_1(X),G)/G \rightarrow H^1(X;G)\]
given by $\wti{F}([\varphi])=[F(\varphi)]$ is the required bijection.
First we need to see that $\wti F([\varphi])$ is well-defined. Take
$g\in G$ and consider the conjugation $g\varphi g^{-1}$. We have
\[F(g\varphi g^{-1})=\gamma\,.\,F(\varphi),\]
where $\gamma\in C^0(X;G)$ is the $0$-cochain which evaluates to~$g$
on each vertex.

Now let us see that $\wti F$ is injective. Assume
$F(\varphi)=\gamma\,.\,F(\psi)$, for some $\gamma\in C^0(X;G)$,
$\varphi,\psi\in\Hom(\pi_1(X),G)$. Since
$F(\varphi)(1,i)=F(\psi)(1,i)=1$, for all $2\leq i\leq n$, we see that
$\gamma$ must have the same value, say $g$, on all the vertices of
$X$. This means that $\varphi=g\psi g^{-1}$, and so
$[\varphi]=[\psi]$.

Finally, let us see that $\wti F$ is surjective. Take an arbitrary
$\sigma\in Z^1(X;G)$.  Define $\gamma\in C^0(X;G)$ by setting
$\gamma(i):=\sigma(1,i)$, for all $2\leq i\leq n$, and
$\gamma(1):=1$. Set $\tau:=\gamma\,.\,\sigma$. Clearly, $\tau(1,i)=1$,
for all $2\leq i\leq n$.  Define $\varphi:\pi_1(X)\ra G$ by setting
$\varphi(e_{ij}):=\tau(i,j)$. Clearly, $F(\varphi)=\tau$, hence $\wti
F([\varphi])=[\sigma]$.  \qed

In particular we obtain the following corollary.

\begin{crl}
\label{hhpp}
Fix an integer $N\geq 2$. The fundamental group $\pi_1(X)$ contains
a~proper normal subgroup $H$, such that $|\pi_1(X):H|\leq N$, if
and only if there exists a~non-trivial simple group $G$, such that
$|G|\leq N$ and the first cohomology $H^1(X;G)$ is non-trivial, i.e.,
$|H^1(X;G)|\geq 2$.
\end{crl}
\pr Assume first that there exists a non-trivial simple group $G$, such
that $|G|\leq N$ and the cohomology group $H^1(X;G)$ is non-trivial.
By Proposition~\ref{prop:fhpi} we know that $|\Hom(\pi_1(X),G)/G|\geq
2$, so we can pick a~non-trivial group homomorphism
$\varphi:\pi_1(X)\ra G$. Set $H:=\ker\varphi$. This is a proper normal subgroup
of $\pi_1(X)$ since $\varphi$ is non-trivial, and $|\pi_1(X):H|=|\im\varphi|\leq
|G|\leq N$.

In the opposite direction, assume that the fundamental group
$\pi_1(X)$ contains a~proper normal subgroup $H$, such that
$|\pi_1(X):H|\leq N$. Let $H$ be chosen so that the index
$|\pi_1(X):H|$ is minimized, and let $G=\pi_1(X)/H$. Then $|G|\leq N$ and $G$ is simple by the choice of $H$. Furthermore, by Proposition~\ref{prop:fhpi}
\[
|H^1(X;G)|=|\Hom(\pi_1(X),G)/G| \geq 2.
\]
 \qed

\mysubs{Non-abelian $1$-Expansion of the Simplex}
\label{expsim}
Let us now adapt our expansion terminology to the non-abelian setting.
Assume $\phi \in C^1(\dn;G)$.  The {\it support} of $\phi$ is
the set
\[\supp\phi:=\Big\{\{u,v\}\in\binom{[n]}{2}:\phi(u,v)\neq 1\Big\}.\]
The {\it norm} of $\phi$ is the cardinality of its support, $\|\phi\|:=|\supp\phi|$.
The {\it cosystolic norm} of $\phi$ is defined as
\[\|\phi\|_{csy}:=\min\{\|\psi\,.\,\phi\|:\psi\in C^0(\dn;G)\}.\]

The following result is an adaptation of Proposition 3.1 of \cite{MW}
to the non-abelian setting.
\begin{prop}
\label{bw1}
Let $\phi \in C^1(\dn;G)$ then
\[\|d_1\phi\| \geq \frac{\ n\, \|\phi\|_{csy}}{3}.\]
\end{prop}
\noindent
{\bf Proof:} For $u \in \dn(0)$ define $\phi_u \in C^0(\dn;G)$ by
setting
\[\phi_u(v):=
\begin{cases}
        \phi(u,v), &\textrm{ if } v \neq u, \\
        1, & \textrm{ otherwise.}
\end{cases}
\]
Note that if $(u,v,w) \in \dn(2)$ then
\begin{equation*}
\begin{split}
(d_1\phi)(u,v,w)&=\phi(u,v)\,\phi(v,w)\,\phi(w,u) \\
&=\phi_u(v)\,\phi(v,w)\,\phi_u(w)^{-1}=(\phi_u\,.\,\phi)(v,w).
\end{split}
\end{equation*}
Therefore
\begin{equation*}
\begin{split}
6\,\|d_1\phi\|&=\big|\{(u,v,w) \in \dn(2):(d_1\phi)(u,v,w) \neq 1\}\big| \\
&=\big|\{(u,v,w) \in \dn(2): (\phi_u\,.\,\phi)(v,w) \neq 1\}\big| \\
&=\sum_{u=1}^n 2\,\|\phi_u\, .\, \phi\| \geq 2n\, \|\phi\|_{csy}.\qed
\end{split}
\end{equation*}

\subsection{Proof of Theorem \ref{nonab}}
\label{fmw}

Let $G$ be an~arbitrary finite group. For a subcomplex
$(\dn)^{(1)} \subset X \subset (\dn)^{(2)}$ we identify
$H^1(X;G)$ with its image under the natural injection $H^1(X;G)
\hookrightarrow H^1\left((\dn)^{(1)};G\right)$. If $\phi \in C^1(\dn;G)$ then
$[\phi] \in H^1(X;G)$ if and only if $(d_1\phi)(u,v,w)=1$ whenever
$(u,v,w) \in\wti X(2)$.  It follows that in the probability space
$Y(n,p)$
\[
\prob\left[ [\phi] \in H^1(Y;G)\right] =(1-p)^{\|d_1\phi\|}.
\]
Therefore, we have
\begin{equation}
\label{unbod}
\begin{split}
\prob \left[H^1(Y;G)\neq\{[1]\}\right]&\leq\sum_{[\phi]}
\prob\left[[\phi] \in H^1(Y;G)\right] \\
&=\sum_{[\phi]}(1-p)^{\|d_1 \phi\|},
\end{split}
\end{equation}
where both sums are taken over all $[\phi]\in H^1\left((\dn)^{(1)};G\right)$,
$[\phi]\neq 1$.

Suppose now that $|G| \leq n^c$. Then
by (\ref{unbod}) and Proposition \ref{bw1} we have
\begin{equation}
\label{unbd1}
\begin{split}
 \prob \left[H^1(Y;G) \neq \{[1]\}\right]  &\leq
\sum_{k \geq 1} \sum_{\|\phi\|_{sys}=k} (1-p)^{\frac{kn}{3}} \\
&\leq \sum_{k \geq 1} \binom{n(n-1)/2}{k} |G|^k
\left(1-\frac{(6+7c)\log n}{n}\right)^{\frac{kn}{3}} \\
&\leq \sum_{k \geq 1} n^{2k} n^{ck} n^{-\frac{(6+7c)k}{3}}=O(n^{-\frac{4c}{3}}).
\end{split}
\end{equation}

Let $\cg(N)$ be the set of all non-trivial simple groups with at most
$N$ elements.  The classification of finite simple groups implies that
there are at most $2$ non-isomorphic simple groups of the same order,
so we certainly have $|\cg(N)| \leq 2N$.  Combining Corollary
\ref{hhpp} and the inequality \eqref{unbd1} we obtain that the
probability that the fundamental group $\pi_1(Y)$ contains a~proper
normal subgroup of index at most $n^c$ cannot exceed
\[
\sum_{G \in \cg(n^c)} \prob\big[Y \in Y(n,p)\,:\,H^1(Y;G) \neq \{[1]\}\big]
\leq O(|\cg(n^c)|n^{-\frac{4c}{3}})=O(n^{-\frac{c}{3}}).\qed
\]

\section{Concluding Remarks}
\label{s:conc}

In this paper we studied several aspects of the $k$-th Cheeger constant of a complex X, a parameter that quantifies the distance of $X$ from a complex $Y$ with nontrivial $k$-th cohomology over $\Int_2$. Our results include, among other things, general methods for bounding the cosystolic norm of a cochain and for bounding the Cheeger constant of a complex, a discussion of expansion of pseudomanifolds and geometric lattices, probabilistic upper bounds on Cheeger constants, and application of non-Abelian expansion to random complexes.
Our work suggests some natural questions regarding higher dimensional expansion:
\begin{itemize}
\item
There are numerous families of combinatorially defined simplicial complexes, e.g.,
chessboard complexes and more general matching complexes, that admit
strong vanishing theorems in (co)homology. It would be interesting to understand whether these vanishing results are accompanied by strong lower bounds on the corresponding Cheeger constants.
\item
In recent years there is a growing interest in developing methodology for studying the topology of
objects (manifolds or more general complexes) using a limited sample of their points. One powerful approach is via persistence homology (see, e.g., Edelsbrunner book \cite{Ed14}). Incorporating Cheeger constants estimates in persistence homology algorithms could lead to improved understanding of the topology of the object. A major challenge in this direction is to devise efficient methods that compute or estimate the expansion of a complex.

\end{itemize}

\end{document}